\documentclass{article}
\usepackage{amssymb}
\usepackage{amsfonts}
\usepackage{amsmath}
\usepackage{amsthm}

\newtheorem{theorem}{Theorem}
\newtheorem{acknowledgement}[theorem]{Acknowledgement}

\newtheorem{claim}[theorem]{Claim}

\newtheorem{conjecture}[theorem]{Conjecture}
\newtheorem{corollary}[theorem]{Corollary}
\newtheorem{criterion}[theorem]{Criterion}

\newtheorem{problem}[theorem]{Problem}
\newtheorem{proposition}[theorem]{Proposition}
\newtheorem{remark}[theorem]{Remark}

\begin{document}

\title{Interpolation and Extrapolation Statements equivalent to the Riemann Hypothesis.}
\author{\'{A}lvaro Corval\'{a}n \\
Instituto del Desarrollo Humano\\
Universidad Nacional de General Sarmiento\\
J. M. Guti\'{e}rrez 1150, C.P. 1613, Malvinas Argentinas\\
Pcia de Bs.As, Rep\'{u}blica Argentina. \\
Mail: acorvala@campus.ungs.edu.ar}
\maketitle

\begin{abstract}
The goal of this paper is twofold; on one hand we wish to present some
statements that can be formulated in terms of Interpolation theory which are
equivalent to the truth or the falseness of the Riemann Hypothesis, on the
other hand we will use a key result of the Jawerth-Milman extrapolation to
improve on a well-known criterion by Beurling and Nyman for the Riemann
Hypothesis giving sharper sufficient conditions to the celebrated hypothesis.

\textit{Keywords: Riemann-Hypothesis, Nyman-Beurling Criterion, Real
Interpolation, Jawerth-Milman Extrapolation, Reverse H\"{o}lder Classes}
\end{abstract}

\section{Introduction}

In 1950 B. Nyman proved that the Riemann hypothesis (RH), which can be
rewritten in terms of the lack of zeros of the Riemann zeta-function in the
half-plane $\operatorname{Re}\left( z\right) >\frac{1}{2}$, is equivalent to the
density in $L^{2}\left( 0,1\right) $ of a certain subspace: 
\begin{equation*} 
\Phi =\left\{ \sum\limits_{k=1}^{n}c_{n_{k}}\left\{ \frac{1}{a_{n_{k}}x}%
\right\} \text{ such that}\sum\limits_{k=1}^{n}c_{n_{k}}\left\{ \frac{1}{%
a_{n_{k}}}\right\} =0\right\} _{a_{n_{k}}\geq 1}\footnote{$\left\{ t\right\}
=t-\left\lfloor t\right\rfloor $ is the fractional part of $t$} 
\end{equation*}%
of this space. A bit later, A. Beurling generalized the result proving that,
in fact, the subspace is free of zeros in the half-plane $\operatorname{Re}%
\left( z\right) >\frac{1}{p}$ if and only if $\Phi $ is dense in $%
L^{p}\left( 0,1\right) $. Since then, a way of trying to prove -or to
disprove- the RH has been to try to demonstrate the density of $\Phi $ in $%
L^{p}$ for all $p<2$ \footnote{%
About fifty years after, L. B\'{a}ez-Duarte strength the Nyman-Beurling
criterion (NB) by showing that it is enough to see the density of the
smaller set, $span\left\{ \left\{ \frac{1}{nx}\right\} \right\} _{n\in 
\mathbb{N}}$ (cf. \cite{BDu5})}. Here we will consider a way of improving
the conditions on the sufficiency side using an extrapolation result to show
that it is enough to prove the density of $\Phi $ in interpolation spaces
between $L^{1}\left( 0,1\right) $ and $L^{2}\left( 0,1\right) $ larger than $%
L^{p}\left( 0,1\right) $ and then with smaller norms.

On the other hand we will show, by means a result of Ivanov and Kalton on
interpolation of subspaces, that if the Riemann hypothesis were not true
this fact can be observed in terms of the Reverse H\"{o}lder class of some
function on $L^{2}\left( 0,1\right) $ through the index characterization of
abstract Reverse H\"{o}lder classes (Cf. \cite{CM}). Specifically we will
see that:

\begin{claim}
If the Riemann Hypothesis is false there is $\psi \in L^{2}\left( 0,1\right) $ 
such that $\psi \in \text{RevH\"{o}lder}\left( L^{2}\left( 0,1\right)
,L^{\infty }\left( 0,1\right) \right) \footnote{%
It is usual to note $RH\left( X,Y\right) $ the class of reverse H\"{o}lder
class for the pair $\left( X,Y\right) $. Nevertheless we will note it
instead $\text{RevH\"{o}lder}(X,Y)$ with the abbreviation RH for Riemann
hypothesis.}$ and $\int_{0}^{1}\psi \left( x\right) \rho \left( \frac{1}{ax}%
\right) dx=0$ for all $a\geq 1$.
\end{claim}

Coming back to the criteria equivalent to RH, one of the ways in which our
main result can be stated is as follows:

\begin{criterion}
\label{Main result} Let $\Phi =\left\{ \sum\limits_{k=1}^{n}c_{n_{k}}\left\{ 
\frac{1}{a_{n_{k}}x}\right\} :\sum\limits_{k=1}^{n}c_{n_{k}}\left\{ \frac{1}{%
a_{n_{k}}}\right\} =0\right\} _{a_{n_{k}}\geq 1}$.

The following statements are equivalent:

a) Riemann Hypothesis is true

b) For $\overline{A}=\left( L^{1}\left( 0,1\right) ,L^{2}\left( 0,1\right)
\right) $ there is a positive function $M\left( \theta \right) $ tempered at 
$1$\footnote{%
We mean that $M\left( \theta \right) $ is tempered at $0$ if $M\left(
2\theta \right) \approx M\left( \theta \right) $ for $\theta $ close to
zero, and $M\left( \theta \right) $ is tempered at $1$ if $M\left( \frac{%
1+\theta }{2}\right) \approx M\left( \theta \right) $ for $\theta $ close to 
$1$. We say that $M$ is tempered if $M$ is tempered at $0$ and at $1$.} with 
$M\left( \theta \right) =O\left( 1\right) $ for $\theta $ going to $1$ such
that%
\begin{equation*}
\inf\limits_{f\in \Phi }\left\Vert \chi _{\left( 0,1\right) }-f\right\Vert
_{\Delta _{0<\theta <1}\left( M\left( \theta \right) \overline{A}_{\theta
,\infty }\right) }=0
\end{equation*}%
Taking $\theta =\frac{2}{p^{\prime }}=2-\frac{2}{p}$ and $\omega \left(
p\right) $ such that $\omega \left( p\right) \left( L^{1},L^{2}\right) _{%
\frac{2}{p^{\prime }},\infty }=M\left( \theta \right) \left(
L^{1},L^{2}\right) _{\theta ,\infty }$ with equivalent norms we can rephrase
b) in the following way:

b') It happens that $\inf\limits_{f\in \Phi }\left\Vert \chi _{\left(
0,1\right) }-f\right\Vert _{\Delta _{1<p<2}\left( \omega \left( p\right)
\left( L^{1},L^{2}\right) _{\frac{2}{p^{\prime }},\infty }\right) }=0$
\end{criterion}

Let's remark that rewriting $K\left( t,f,L^{1},L^{2}\right) $ in terms of $%
K\left( t^{2},f,L^{1},L^{\infty }\right) $ and using that in this context it
is possible to consider the modified spaces $\left\langle \overline{A}%
_{\theta ,q,K}\right\rangle $ (see the definition in the next section)
instead of the interpolation spaces $\overline{A}_{\theta ,q,K}$, then the
conditions b) and b') respectively mean that: 
\begin{equation} 
\inf\limits_{f\in \Phi }\left( \sup\limits_{0<\theta <1}\left[
\sup\limits_{0<t<1}M\left( \theta \right) t^{1-\theta }\left(
\int\limits_{t^{2}}^{\infty }\left( \chi _{\left( 0,1\right) }-f\right)
^{\ast \ast }\left( s\right) ^{2}ds\right) ^{\frac{1}{2}}\right] \right) =0
\label{D1} 
\end{equation}%
and 
\begin{equation} 
\inf\limits_{f\in \Phi }\left( \sup\limits_{1<p<2}\left[ \sup\limits_{0<t<1}%
\omega \left( p\right) t^{1-\frac{2}{p^{\prime }}}\left(
\int\limits_{t^{2}}^{\infty }\left( \chi _{\left( 0,1\right) }-f\right)
^{\ast \ast }\left( s\right) ^{2}ds\right) ^{\frac{1}{2}}\right] \right) =0
\label{D2} 
\end{equation}%
, and then both \ref{D1} and \ref{D2} are equivalent to the Riemann
Hypothesis. The details will appear below.

\section{Preliminaries and description by of the problem by means of
Interpolation}

Now we wish to set a kind of problems in terms of interpolation spaces.

A couple of Banach spaces $X_{0}$ and $X_{1}$ is a compatible pair $%
\overline{X}=\left( X_{0},X_{1}\right) $ if there is a Hausdorff topological
vector space $\mathcal{H}$ with $X_{0},X_{1}\subset \mathcal{H}$ with
continuous inclusions. A morphism $T:\overline{X}\rightarrow \overline{Y}$
between compatible pairs of Banach spaces $\overline{X}=\left(
X_{0},X_{1}\right) $ and $\overline{Y}=\left( Y_{0},Y_{1}\right) $ is a
bounded linear operator $T$ from $X_{0}+X_{1}$ into $Y_{0}+Y_{1}$ with
bounded restrictions of $T$ mapping $X_{i}$ into $Y_{i}$ for $i=0,1$.

Given $\overline{A}=\left( A_{0},A_{1}\right) $ and $\overline{B}=\left(
B_{0},B_{1}\right) $ compatible pairs of Banach spaces if $A_{0}\cap
A_{1}\subset A\subset A_{0}+A_{1}$ and $B_{0}\cap B_{1}\subset B\subset
B_{0}+B_{1}$ in both cases with continuous embeddings, we say that $A$ and $%
B $ are interpolation spaces relative to $\overline{A}$ and $\overline{B}$
if for every $T:\overline{A}\rightarrow \overline{B}$ has an extension $%
T:A\rightarrow B$. If $\left\Vert T\right\Vert _{A\rightarrow B}\leq \max
\left( \left\Vert T\right\Vert _{A_{0}\rightarrow B_{0}},\left\Vert
T\right\Vert _{A_{1}\rightarrow B_{1}}\right) $ we say that $A,B$ are exact.

An interpolation method $\mathcal{F}$ is a functor from the category of
compatible pairs of Banach spaces into the category of Banach spaces such
that if $\overline{A}=\left( A_{0},A_{1}\right) $ and $\overline{B}=\left(
B_{0},B_{1}\right) $ are pairs of Banach spaces then $\mathcal{F}\left( 
\overline{A}\right) $ and $\mathcal{F}\left( \overline{B}\right) $ are
interpolation spaces relative to $\overline{A}$ and $\overline{B}$, that is\ 
$A_{0}\cap A_{1}\subset \mathcal{F}\left( \overline{A}\right) \subset
A_{0}+A_{1}$, $B_{0}\cap B_{1}\subset \mathcal{F}_{\theta }\left( \overline{B%
}\right) \subset B_{0}+B_{1}$ and we have $T:\mathcal{F}\left( \overline{A}%
\right) \rightarrow \mathcal{F}\left( \overline{B}\right) $ whenever $T:%
\overline{A}\rightarrow \overline{B}$. If $\mathcal{F}$ always yields exact
interpolation spaces $\mathcal{F}\left( \overline{A}\right) ,\mathcal{F}%
\left( \overline{B}\right) $ we say that $\mathcal{F}$ is exact. If,
moreover, $\left\Vert T\right\Vert _{\mathcal{F}\left( A\right) \rightarrow 
\mathcal{F}\left( B\right) }\leq \left( \left\Vert T\right\Vert
_{A_{0}\rightarrow B_{0}}\right) ^{1-\theta }\left( \left\Vert T\right\Vert
_{A_{1}\rightarrow B_{1}}\right) ^{\theta }$ for some fixed $\theta \in
\left( 0,1\right) $ then $\mathcal{F}$ is called exact of exponent $\theta $.

Let $\left( A_{0},A_{1}\right) $ an ordered pair of Banach spaces, that is $%
A_{1}\overset{1}{\subset }A_{0}\footnote{$A_{1}\overset{1}{\subset }A_{0}$
means that the norm of the inclusion is $\leq 1$. It's immediate that for ordered pairs we don't need to ask that both $%
A_{0},A_{1}$ to be included in some Hausdorff vector space $\mathcal{H}$ and
it is clear that $A_{0}+A_{1}=A_{0}$ and $A_{0}\cap A_{1}=A_{1}$}$ and let $%
\left\{ f_{i}\right\} _{i\in I}\subset A_{1}$ a family of functions dense in 
$A_{0}$. We wish to give conditions that ensure that $\Phi =span\left\{
f_{i}\right\} _{i\in I}$ is also dense in $A_{1}\footnote{%
We ask for the density of $\Phi $ with the norm of $A_{1}$. Being $\Phi $
dense in $A_{0}$, for $f\in A_{1}\subset A_{0}$ and any $\varepsilon >0$ we
have $f_{\varepsilon }\in \Phi $ such that $\left\Vert f-f_{\varepsilon
}\right\Vert _{A_{0}}<\varepsilon $ but we don't know if we can find some $\widetilde{f}_{\varepsilon }\in \Phi $ such that 
$\left\Vert f-\widetilde{f}_{\varepsilon }\right\Vert _{A_{1}}<\varepsilon $.%
}$. As it is obvious, a necessary condition is that $A_{1}$ itself and a
fortiori $A_{0}\cap A_{1}$ must be dense in $A_{0}$, so we will assume that $%
\left( A_{0},A_{1}\right) $ is a regular pair, which it means that $%
A_{0}\cap A_{1}$ is dense in both $A_{0}$ and $A_{1}$.

An ordered scale is a family $\left\{ A_{\theta }\right\} _{\theta \in
\Theta }\footnote{%
Usually the index set is $\Theta =\left( 0,1\right) $}$ of Banach spaces
such that $A_{\theta _{1}}\subset A_{\theta _{2}}$ whenever $\theta
_{2}<\theta _{1}$. We will mostly deal with ordered scales given by a family
of interpolation functors $\left\{ \mathcal{F}_{\theta }\right\} _{\theta
\in \Theta }$ applied on a compatible ordered pair $\left(
A_{0},A_{1}\right) $ with $A_{1}\subset A_{0}$ such that $\left\{ \mathcal{F}%
_{\theta }\left( A\right) \right\} _{\theta \in \Theta }$ is an ordered
scale, that is $A_{1}\overset{1}{\subset }\mathcal{F}_{\theta }\left( 
\overline{A}\right) \overset{1}{\subset }A_{0}$ and $\mathcal{F}_{\theta
_{2}}\left( \overline{A}\right) \overset{1}{\subset }\mathcal{F}_{\theta
_{1}}\left( \overline{A}\right) $ for $0\leq \theta _{1}<\theta _{2}\leq 1$
we would also be interested in knowing how to characterize the set $C_{%
\mathcal{F}}=\{\theta :\Phi $ is dense in $\mathcal{F}_{\theta }\left( 
\overline{A}\right) \}$. In other words, we would like to find out when $%
\overline{\Phi }^{\mathcal{F}_{\theta }\left( \overline{A}\right) }=\mathcal{%
F}_{\theta }\left( \overline{A}\right) $ where $\overline{\Phi }^{\mathcal{F}%
_{\theta }\left( \overline{A}\right) }$ means the closure of $\Phi $ in $%
\mathcal{F}_{\theta }\left( \overline{A}\right) $.

\section{Extrapolation and improvements}

We will show that we can use, in a straightforward way, a key fact of the
Jawerth-Milman extrapolation theory to improve a density criterion for a
family of interpolation spaces, one of which is of particular importance. We
begin by describing such instance:

First we remember that the famous Riemann Hypothesis (RH) conjectures that
the Riemann zeta-function (that is the analytic continuation of the function
defined by $\zeta \left( s\right) =\sum\limits_{n=1}^{\infty }\frac{1}{n^{s}}$
for $\operatorname{Re}\left( s\right) >1$) has non trivial zeros only in the
critical line $\operatorname{Re}\left( s\right) =\frac{1}{2}$ where the trivial
zeros are the negative even integers. It is known that any non-trivial zero
of $\zeta \left( s\right) $ lies in the open strip$\ \left\{ s\in \mathbb{C}%
:0<\operatorname{Re}\left( s\right) <1\right\} $ known as the critical strip, and
the functional identity 
\begin{equation*} 
\zeta \left( s\right) =2^{s}\pi ^{s-1}\sin \left( \frac{\pi s}{2}\right)
\Gamma \left( 1-s\right) \zeta \left( 1-s\right) 
\end{equation*}%
for any $s\in \mathbb{C}$ except for $s=0$ or $1$ shows us that the
non-trivial zeros of $\zeta $ are symmetric respect to the critical line $%
\operatorname{Re}\left( s\right) =\frac{1}{2}$ and then we only have to check the
existence or or lackness of non-trivial zeros of $\zeta $ in $\left\{ \frac{1%
}{2}<\operatorname{Re}\left( s\right) <1\right\} $ because any hypothetical
non-trivial zero in $\left\{ 0<\operatorname{Re}\left( s\right) <\frac{1}{2}\right\} 
$ will be paired with another in $\left\{ \frac{1}{2}<\operatorname{Re}\left(
s\right) <1\right\} $. Also it is well known that $\zeta $ do have infinite
non-trivial zeros in the critical line $\left\{ \operatorname{Re}\left( s\right) =%
\frac{1}{2}\right\} $. Thus the Riemann Hypothesis can be stated as follows
(cf. \cite{R}):

\begin{conjecture}
(Bernhard Riemann) \textit{Riemann zeta-function }$\zeta $\textit{\ is free
from zeros for }$\frac{1}{2}<\operatorname{Re}\left( s\right) <1$.
\end{conjecture}

Secondly, there are well known results about $\zeta $ that allows to state
RH in terms of a density problem (cf. \cite{B},\cite{N}, and \cite{BDu}):

\begin{theorem}
\label{Nyman-Beurling}(Nyman-Beurling) Let $\rho :\left( 0,1\right)
\rightarrow \mathbb{R}_{\geq 0}$, the fractional part of $x$, given by $\rho
\left( x\right) =x-\left\lfloor x\right\rfloor $ and $\rho _{a}\left(
x\right) =\rho \left( \frac{1}{ax}\right) $ for $1\leq a<\infty $. The
Riemann zeta-function $\zeta $ is zero-free in the half-plane $\sigma >%
\frac{1}{p}$ where $s=\sigma +i\tau $ for $1\leq p<\infty $, if and only if $%
\Phi =\left\{ \sum\limits_{k=1}^{n}c_{n_{k}}\rho \left( \frac{1}{a_{n_{k}}x}%
\right) :\sum\limits_{k=1}^{n}c_{n_{k}}\rho \left( \frac{1}{a_{n_{k}}}%
\right) =0\right\} _{a_{n_{k}}\geq 1}$ is dense in the space $L^{p}\left(
0,1\right) $.
\end{theorem}

\begin{remark}
As we have already mentioned, Luis B\'{a}ez-Duarte proved that it is enough
to consider only integer positive dilations values of $a$, that is
considering 
\begin{equation*} 
\Phi =\left\{ \sum\limits_{k=1}^{n}c_{n_{k}}\rho \left( \frac{1}{n_{k}x}%
\right) :\sum\limits_{k=1}^{n}c_{n_{k}}\rho \left( \frac{1}{n_{k}}\right)
=0\right\} _{n_{k}\in \mathbb{N}} 
\end{equation*}
instead of 
\begin{equation*} 
\left\{ \sum\limits_{k=1}^{n}c_{n_{k}}\rho \left( \frac{1}{a_{n_{k}}x}%
\right) :\sum\limits_{k=1}^{n}c_{n_{k}}\rho \left( \frac{1}{a_{n_{k}}}%
\right) =0\right\} _{a_{n_{k}}\geq 1} 
\end{equation*}%
. In any case our argument to improve the criterion works, so in the
following we will only write $\Phi $.
\end{remark}

\begin{remark}
The interesting problem about which functions have the property that the
span of their dilations generate $L^{2}$ was raised in the fifties by
Wintner, Beurling, Kozlov and others and studied up today (cf. \cite{K}, 
\cite{HK}, \cite{Y}, and \cite{Mit})
\end{remark}

Although it seems, at first glance, that for proving RH it is necessary to
check out if that for every function $f$ in $L^{p}$ we can find $%
f_{n}\rightarrow f$ on $L^{p}$ with $f_{n}=\sum\limits_{k=1}^{n}c_{n_{k}}%
\rho \left( \frac{1}{a_{n_{k}}x}\right) $, actually we only need to verify
the statement for $f=\chi _{\left( 0,1\right) }$, as it is noted in
Beurling's celebrated paper (cf. \cite{B}), because if $\chi _{\left(
0,1\right) }\in \overline{\Phi }$ so it is $\chi _{\left( 0,a\right) }\in 
\overline{\Phi }$ for $a>0$, and so they are the characteristics of
intervals included in $\left( 0,1\right) $, and then any step function and
therefore, by density, the whole $L^{p}$. That is we have the following:

\begin{proposition}
$\Phi $ is dense in $L^{p}\left( 0,1\right) $ if and only if $\overline{\Phi 
}^{L^{p}}$(the closure of $\Phi $ in $L^{p}$) contains $\chi _{\left(
0,1\right) }$.
\end{proposition}

Thus we can use the Nyman-Beurling Theorem (\ref{Nyman-Beurling}) to obtain
the following equivalence (cf. \cite{B},\cite{N},\cite{BDu4}):

\begin{theorem}
\label{N-B equivalence}The following statements are equivalent:

i) \textit{Riemann zeta-function }$\zeta $\textit{\ is zero-free for }$%
\frac{1}{2}<\operatorname{Re}\left( s\right) <1$, or equivalently the Riemann
Hypothesis is true.

ii) $\chi _{\left( 0,1\right) }\in \overline{\Phi }^{L^{2}}$.

iii) $\chi _{\left( 0,1\right) }\in \overline{\Phi }^{L^{p}}$ for every $%
p\in (1,2)$.$\footnote{%
Taking in account that $\Phi $ is close by dilations with scale factor in $%
(0,1]$ it is clear that it is possible to consider any other generator
instead of $\chi \left( t\right) =\chi _{\left( 0,1\right) }$, where we mean
by generator a function in $L^{p}\left( 0,1\right) $ which span of its
dilations with scale factor in $[1,\infty )$ is dense in $L^{p}$. For
instance $\lambda \left( t\right) =\chi \left( t\right) \log \left( t\right) 
$ is also a well known generator (cf. \cite{BDu4}).}$
\end{theorem}

\begin{remark}
The equivalence between i) and iii) and between i) and ii) is guaranteed by
the Nyman-Beurling Theorem \ref{Nyman-Beurling} and the openness of the
interval $\left( \frac{1}{2},1\right) $.
\end{remark}

\begin{corollary}
We can write the equivalence of the Riemann Hypothesis with the statements
ii) and iii) in terms of the distances in $L^{p}$ from $\chi _{\left(
0,1\right) }$ to the closure of $\Phi $.
\end{corollary}

\begin{corollary}
\label{NBC}A) The Riemann Hypothesis is true if and only if 
\begin{equation*} 
\inf\limits_{f_{n}\in \Phi }\left\Vert \chi _{\left( 0,1\right)
}-f_{n}\right\Vert _{p}=0 
\end{equation*}%
for every $p\in \left( 1,2\right) $ , that is: 
\begin{equation*} 
RH\Longleftrightarrow \lim\limits_{n\rightarrow \infty }\left(
\inf\limits_{a_{n_{k}},c_{n_{k}}}\left( \int_{0}^{1}\left\vert
1-\sum\limits_{k=1}^{n}c_{n_{k}}\rho \left( \frac{1}{a_{n_{k}}x}\right)
\right\vert ^{p}dx\right) ^{\frac{1}{p}}\right) =0 
\end{equation*}%
for any $p\in \left( 1,2\right) $ where the infimum is taken on $%
a_{n_{k}}\geq 1$ and $c_{n,k}\in \mathbb{R}$.

B) The Riemann Hypothesis is true if and only if 
\begin{equation*} 
\inf\limits_{f_{n}\in \Phi }\left\Vert \chi _{\left( 0,1\right)
}-f_{n}\right\Vert _{2}=0 
\end{equation*}%
, that is: 
\begin{equation*} 
RH\Longleftrightarrow \lim\limits_{n\rightarrow \infty }\left(
\inf\limits_{a_{n_{k}},c_{n_{k}}}\left( \int_{0}^{1}\left\vert
1-\sum\limits_{k=1}^{n}c_{n_{k}}\rho \left( \frac{1}{a_{n_{k}}x}\right)
\right\vert ^{2}dx\right) ^{\frac{1}{2}}\right) =0 
\end{equation*}%
where the infimum is taken on $a_{n_{k}}\geq 1$ and $c_{n,k}\in \mathbb{R}$.
\end{corollary}

\subsection{A bit of Extrapolation Theory}

In Extrapolation Theory we deal with strongly compatible families\footnote{%
From now on we will call them compatible families.} of Banach spaces indexed
by some fixed index set $\left\{ A_{\theta }\right\} _{\theta \in \mathbf{%
\Theta }}$, usually $\Theta =\left( 0,1\right) $. We will say that a family $%
A=\left\{ A_{\theta }\right\} _{\theta \in \Theta }$ of Banach spaces is
compatible if there exists two Banach spaces $\mathbb{A}_{0}$ and $\mathbb{A}%
_{1}$ such that $\mathbb{A}_{1}\subset A_{\theta }\subset \mathbb{A}_{0}$
with continuous inclusions for each $\theta \in \Theta $. For a fixed index
set $\Theta $ we will consider the category of compatible families of Banach
spaces whose morphisms are the linear operators $T:\left\{ A_{\theta
}\right\} _{\theta \in \Theta }\rightarrow \left\{ B_{\theta }\right\}
_{\theta \in \Theta }$ which means that there is a linear operator $T:%
\mathbb{A}_{0}\rightarrow \mathbb{B}_{0}$ with restrictions $T:A_{\theta }%
\overset{1}{\rightarrow }B_{\theta }\footnote{%
This means that $\left\Vert T\right\Vert _{A_{\theta }\rightarrow B_{\theta
}}\leq 1$}$ for every $\theta \in \Theta $. We say that $A$ and $B$ are
extrapolation spaces respect to $\left\{ A_{\theta }\right\} _{\theta \in
\Theta }$ and $\left\{ B_{\theta }\right\} _{\theta \in \Theta }$ if $%
\mathbb{A}_{1}\subset A\subset \mathbb{A}_{0}$, $\mathbb{B}_{1}\subset
B\subset \mathbb{B}_{0}$ and $T:A\rightarrow B$ whenever $T:\left\{
A_{\theta }\right\} _{\theta \in \Theta }\rightarrow \left\{ B_{\theta
}\right\} _{\theta \in \Theta }$. An extrapolation method $\mathcal{E}$ is a
functor from a collection of compatible families $Dom\left( \mathcal{E}%
\right) $ into a collection of Banach spaces such that $\mathcal{E}\left(
\left\{ A_{\theta }\right\} _{\theta \in \Theta }\right) $ and $\mathcal{E}%
\left( \left\{ B_{\theta }\right\} _{\theta \in \Theta }\right) $ are
extrapolation spaces respect to $\left\{ A_{\theta }\right\} _{\theta \in
\Theta }$ and $\left\{ B_{\theta }\right\} _{\theta \in \Theta }$. The
simplest but also more important extrapolation functors are the $\Sigma $
and $\Delta $ methods:

For the compatible families such that the injections $\mathbb{A}%
_{1}\hookrightarrow A_{\theta }$ are uniformly bounded, that is $%
\sup\limits_{\theta \in \Theta }\sup\limits_{a\in \mathbb{A}_{1}}\frac{%
\left\Vert a\right\Vert _{A_{\theta }}}{\left\Vert a\right\Vert _{\mathbb{A}%
_{1}}}<\infty $, we define 
\begin{equation*} 
\Delta \left( \left\{ A_{\theta }\right\} _{\theta \in \Theta }\right)
=\left\{ a\in \bigcap\limits_{\theta \in \Theta }A_{\theta }:\left\Vert
a\right\Vert _{\Delta \left\{ A_{\theta }\right\} _{\theta \in \Theta
}}:=\sup\limits_{\theta \in \Theta }\left\Vert a\right\Vert _{A_{\theta
}}<\infty \right\} 
\end{equation*}%
, and then $\Delta $ provides us a functor from this category to the
category of Banach spaces (cf. \cite{JM}).

The dual $\Sigma $ method is defined for the compatible families such that
the injections $A_{\theta }\hookrightarrow \mathbb{A}_{1}$ such that $%
\sup\limits_{\theta \in \Theta }\sup\limits_{a\in A_{\theta }}\frac{%
\left\Vert a\right\Vert _{\mathbb{A}_{0}}}{\left\Vert a\right\Vert
_{A_{\theta }}}<\infty $, in this context 
\begin{equation*} 
\Sigma \left( \left\{ A_{\theta }\right\} _{\theta \in \Theta }\right)
=\left\{ a\in \mathbb{A}_{0}:\left\Vert a\right\Vert _{\Sigma \left\{
A_{\theta }\right\} _{\theta \in \Theta }}:=\inf \sum\limits_{\theta
}\left\Vert a_{\theta }\right\Vert _{A_{\theta }}<\infty \right\} 
\end{equation*}%
where the infimum is taken over all representations $a=\sum a_{\theta
}:a_{\theta }\in A_{\theta }$ such that $\sum\limits_{\theta }\left\Vert
a_{\theta }\right\Vert _{A_{\theta }}<\infty $, understanding, as usual that 
$\inf \emptyset =\infty $.

Under quite general conditions we have a close connection between
extrapolation and interpolation: Let's consider compatible families $\left\{ X_{\theta }\right\} _{\theta \in
\Theta }$ with $X_{\theta }=\mathcal{F}_{\theta }\left( \overline{X}\right) $
for a family of interpolation functors $\left\{ \mathcal{F}_{\theta
}\right\} _{\theta \in \Theta }$ applied over a compatible pair of Banach
spaces $\overline{X}=\left( X_{0},X_{1}\right) $. There is a
characterization due to Jawerth and Milman (cf. \cite{JM}) of the families
of methods $\left\{ \mathcal{F}_{\theta }\right\} _{\theta \in \Theta }$
such that we can "reverse" in certain sense the extrapolation property: A
family $\left\{ \mathcal{F}_{\theta }\right\} _{\theta \in \Theta }$ is
complete if $T:\mathcal{F}_{\theta }\left( \overline{A}\right) \overset{1}{%
\rightarrow }\mathcal{F}_{\theta }\left( \overline{B}\right) $ for all $%
\theta \in \Theta $ implies that $T:\overline{A}\rightarrow \overline{B}$
whenever $\overline{A}=\left( A_{0},A_{1}\right) $ is a regular pair and $%
\overline{A}$ and $\overline{B}$ are mutually closed\footnote{$\overline{X}%
=\left( X_{0},X_{1}\right) $ is mutually close if $X_{0}=\widetilde{X_{0}}$
where $\widetilde{X_{0}}=\left\{ f\in X_{0}+X_{1}:\lim\limits_{t\rightarrow
\infty }\left( \inf\limits_{f_{0},f_{1}:f=f_{0}+f_{1}}\left[ \left\Vert
f_{0}\right\Vert _{X_{0}}+t\left\Vert f_{1}\right\Vert _{X_{1}}\right]
\right) <\infty \right\} $ and $X_{1}=\widetilde{X_{1}}$ (defined
analogously).}. For to describe complete families, characteristic functions
are defined by the action of the functors $\mathcal{F}$ over one dimensional
spaces: the characteristic function of $\mathcal{F}$ is $\rho $ given by $%
\mathcal{F}\left( \left( \mathbb{C},\frac{1}{t}\mathbb{C}\right) \right) =%
\frac{1}{\rho \left( t\right) }\mathbb{C},t>0$. It is straightforward that
for $\mathcal{F}$ exact the function $\rho $ is quasi-concave\footnote{%
This means that $\rho _{\theta }$ is increasing and $\frac{\rho _{\theta
}\left( t\right) }{t}$ is decreasing as function of $t$)}. Moreover, if $%
\mathcal{F}$ is exact of exponent $\theta $ then $\rho \left( t\right)
=Ct^{\theta }$ for some $C>0$. Now, the characterization mentioned above is
the following theorem for which proof we refer the reader to \cite{JM}:

\begin{theorem}
(\cite{JM} Th. 2.5) A family of exact interpolation functors $\left\{ 
\mathcal{F}_{\theta }\right\} _{\theta \in \Theta }$ is complete if and only
if $\inf\limits_{\theta \in \Theta }\frac{\rho _{\theta \left( t\right) }}{%
\rho _{\theta \left( s\right) }}\leq C\min \left( 1,\frac{t}{s}\right) $ for
all $s,t>0$.
\end{theorem}

\begin{remark}
From the above theorem it is immediate that if each $\mathcal{F}_{\theta }$
is exact of exponent $\theta $ then $\left\{ \mathcal{F}_{\theta }\right\}
_{\theta \in \Theta }$ is complete.
\end{remark}

\begin{remark}
Among the exact interpolation methods with a fixed characteristic function $%
\rho $ there are well known the minimal and maximal methods. Let's describe this situation in terms of the so called real methods of
interpolation:

For a compatible pair of Banach spaces $\overline{X}=\left(
X_{0},X_{1}\right) $ and $f\in X_{0}+X_{1}$, the $K-functional$ is given by 
\begin{equation*} 
K\left( t,f,X_{0},X_{1}\right) =\inf\limits_{f_{0}+f_{1}=f}\left\{
\left\Vert f_{0}\right\Vert +t\left\Vert f_{1}\right\Vert :f_{0}\in
X_{0},f_{1}\in X_{1}\right\} 
\end{equation*}%
for $t>0$. On the other hand for $f\in X_{0}\cap X_{1}$ the $J-functional$
is given by $J\left( t,f,X_{0},X_{1}\right) =\max \left\{ \left\Vert
f\right\Vert _{X_{0}},t\left\Vert f\right\Vert _{X_{1}}\right\} $ for any $%
t>0$. Now, if $\rho $ is quasi-concave and $1\leq q\leq \infty $, we define
the interpolation space 
\begin{equation*} 
\overline{X}_{\rho ,q,K}=\left( X_{0},X_{1}\right) _{\rho ,q,K}=\left\{ f\in
X_{0}+X_{1}:\left\Vert f\right\Vert _{\rho ,q,K}<\infty \right\} 
\end{equation*}%
where \ $\left\Vert f\right\Vert _{\rho ,q,K}:=\phi _{\theta ,q}\left(
K\left( t,f,X_{0},X_{1}\right) \right) $, being%
\begin{equation*} 
\phi _{\theta ,q}\left( \psi \right) =\left\{ 
\begin{array}{c}
\left( \int_{0}^{\infty }\left( \frac{\psi \left( s\right) }{\rho \left(
s\right) }\right) ^{q}\frac{ds}{s}\right) ^{\frac{1}{q}}\text{ if }q<\infty
\\ 
\sup\limits_{s>0}\frac{\psi \left( s\right) }{\rho \left( s\right) }\text{
if }q=\infty%
\end{array}%
\right. 
\end{equation*}%
. Whereas for the $J-functional$ we have the interpolation spaces: 
\begin{equation*} 
\overline{X}_{\rho ,q,J}=\left( X_{0},X_{1}\right) _{\rho ,q,J}=\left\{ f\in
X_{0}+X_{1}:\left\Vert f\right\Vert _{\rho ,q,J}<\infty \right\} 
\end{equation*}%
where $\left\Vert f\right\Vert _{\rho ,q,J}=\inf\limits_{u}\left\{ \phi
_{\theta ,q}\left( J\left( t,u\left( t\right) ,X_{0},X_{1}\right) \right)
\right\} $ where the infimum is taken over the representations $u:\left(
0,\infty \right) \rightarrow X_{0}\cap X_{1}$ such that $\int_{0}^{\infty
}u\left( t\right) \frac{dt}{t}=f$ with convergence in $X_{0}+X_{1}$.

Note that $\overline{X}_{\rho ,q,K}$ is an intermediate space for all $\rho $
quasiconcave as long as $q=\infty $ and so it is $\overline{X}_{\rho ,q,J}$
provided that $q=1$.

At last, we can tell about the minimal and the maximal methods for $\rho $:
\end{remark}

\begin{theorem}
If $\mathcal{F}$ is an exact interpolation method with characteristic
function $\rho $ and $\overline{X}$ is a compatible pair of Banach spaces,
then 
\begin{equation*} 
\overline{X}_{\rho ,1,J}\overset{1}{\rightarrow }\mathcal{F}\left( \overline{%
X}\right) \overset{1}{\rightarrow }\overline{X}_{\rho ,\infty ,K} 
\end{equation*}
\end{theorem}

We refer, for instance, to \cite{AM} for the proof of the above theorem
which is essentially due to Lions and Peetre -for $\rho \left( t\right)
=Ct^{\theta },0<\theta <1$ and S. Janson for the slightly improvement with
any $\rho $ quasi-concave-.

\begin{remark}
When $\rho \left( t\right) =Ct^{\theta }$ we get the notation $\overline{X}%
_{\theta ,q,K}$ and $\overline{X}_{\theta ,q,J}$ is adopted instead of $%
\overline{X}_{\rho ,q,K}$ and $\overline{X}_{\rho ,q,J}$. In Extrapolation
Theory, it's usual to modify the norm of $\overline{X}_{\theta ,q,K}$ and $%
\overline{X}_{\theta ,q,J}$, multiplying by suitable constants to ensure
that the norm is exactly $\rho \left( t\right) =t^{\theta }$. We define $%
\left\Vert f\right\Vert _{\theta ,q,K}^{\blacktriangleleft }=\left( \left(
1-\theta \right) \theta q\right) ^{\frac{1}{q}}\left\Vert f\right\Vert
_{\theta ,q,K}$ with the convention: $\left( \left( 1-\theta \right) \theta
q\right) ^{\frac{1}{q}}=1$ for $q=\infty $ and $\left\Vert f\right\Vert
_{\theta ,q,J}^{\blacktriangleleft }=\left( \left( 1-\theta \right) \theta
q^{\prime }\right) ^{-\frac{1}{q^{\prime }}}\left\Vert f\right\Vert _{\theta
,q,J}$ with $\left( \left( 1-\theta \right) \theta q^{\prime }\right) ^{-%
\frac{1}{q^{\prime }}}=1$ for $q=1$. We also adopt the notation $\overline{X}%
_{\theta ,q,K}^{\blacktriangleleft }$ and $\overline{X}_{\theta
,q,J}^{\blacktriangleleft }$ for the spaces obtained from $\overline{X}%
_{\theta ,q,K}$ and $\overline{X}_{\theta ,q,J}$ respectively renormed with $%
\left\Vert f\right\Vert _{\theta ,q,K}^{\blacktriangleleft }$ or $\left\Vert
f\right\Vert _{\theta ,q,J}^{\blacktriangleleft }$. Let's observe that with the above convention $\overline{X}_{\theta ,\infty
,K}^{\blacktriangleleft }=\overline{X}_{\theta ,\infty ,K}$, so in this case
we keep the notation $\overline{X}_{\theta ,\infty ,K}$.
\end{remark}

\begin{remark}
The equivalence between ii) and iii) in theorem \ref{N-B equivalence} for $%
L^{p}\left( 0,1\right) $ can be consider in terms of the extrapolation
functor-$\Delta $. Let's explain this statement: First of all, $L^{p}\left( 0,1\right) $ spaces can be
obtained as interpolation spaces of the pair $\left( L^{1}\left(
0,1\right) ,L^{2}\left( 0,1\right) \right) $. To see that, we first recall
some well-known facts about the pair $\left( L^{1},L^{\infty }\right) $: On
one hand we can obtain the Lorentz spaces $L^{p,q}$ as interpolation spaces
for the pair $\left( L^{1},L^{\infty }\right) _{\theta ,q}\footnote{%
As is well known, for $0<\theta <1$ and $1\leq q\leq \infty $ we have that $%
\overline{X}_{\theta ,q,K}=\overline{X}_{\theta ,q,J}$ with equivalence of
norms. Hence, in such cases, we omit the reference to the interpolation
method.}$ for $p=\frac{1}{1-\theta }$ (cf. \cite{BS}, pag. 300) and, in
particular, $L^{p}=L^{p,p}=\left( L^{1},L^{\infty }\right) _{\frac{1}{%
p^{\prime }},p}$, on the other hand, the interpolation spaces of the pair $%
\left( L^{1},L^{2}\right) $ can be obtained by using reiteration from the
pair $\left( L^{1},L^{\infty }\right) $ as stated in \cite{BS} (theorem
V.2.4, pag 311): For this we name $\overline{X}=\left( X_{0},X_{1}\right) $
for $X_{0}=L^{1}$, $X_{1}=L^{\infty }$, and we consider the intermediate
spaces: $\overline{X}_{\theta _{0}}=L^{1}=\left( L^{1},L^{\infty }\right)
_{0,1}$ which is of class $\theta _{0}=0\footnote{%
An intermediate space $X$ of a compatible pair $\left( X_{0},X_{1}\right) $
is said to be of class $\theta $ with $0<\theta <1$ if $\left(
X_{0},X_{1}\right) _{\theta ,1}\hookrightarrow X\hookrightarrow \left(
X_{0},X_{1}\right) _{\theta ,\infty }$. It is said of class $\theta =0$ if $%
X_{0}\hookrightarrow X\hookrightarrow \widetilde{X_{0}}$, and it is said of class $\theta =1$ if $X_{1}\hookrightarrow
X\hookrightarrow \widetilde{X_{1}}$.}$ and $\overline{X}_{\frac{1}{2}%
}=L^{2}=\left( L^{1},L^{\infty }\right) _{\frac{1}{2},2}$ which is of class $%
\theta _{1}=\frac{1}{2}$, so we have that $\left( \overline{X}_{\theta _{0}},%
\overline{X}_{\theta _{1}}\right) _{\theta ,p}=\left( X_{0},X_{1}\right)
_{\eta ,p}$ for $\eta =\left( 1-\theta \right) \theta _{0}+\theta \theta
_{1} $ which for the pair $\left( L^{1},L^{\infty }\right) $ and for $\theta
_{0}=0$ and $\theta _{1}=\frac{1}{2}$ means that $\left( L^{1},L^{2}\right)
_{\theta ,p}=\left( L^{1},L^{\infty }\right) _{\eta ,p}$ for $\eta =\frac{1}{%
2}\theta $. As a consequence of this we have that $L^{p}=\left(
L^{1},L^{\infty }\right) _{\frac{1}{p^{\prime }},p}$ is also an
interpolation space of $\left( L^{1},L^{2}\right) $, specifically $\left(
L^{1},L^{2}\right) _{\theta ,\frac{2}{2-\theta }}$ for $p=\frac{2}{2-\theta }
$ and $\theta =\frac{2}{p^{\prime }}=2\left( 1-\frac{1}{p}\right) $.Then, $%
\left( L^{1},L^{2}\right) _{\theta ,\frac{2}{2-\theta }}$ runs along $L^{p}$
spaces with $p\in \left( 1,2\right) $ while $\theta $ goes from $0$ to $1$.

Finally, for $p=\frac{2}{2-\theta }$, we renorm $\left( L^{1},L^{2}\right)
_{\theta ,p}=L^{p}\left( 0,1\right) $ spaces with norm $\left( \left(
1-\theta \right) \theta \frac{2}{2-\theta }\right) ^{\frac{2-\theta }{2}%
}\left\Vert f\right\Vert _{\theta ,\frac{2}{2-\theta },K}$ as in the
previous remark, because we wish to use a family $\mathcal{F}_{\theta
}\left( \overline{A}\right) $ of interpolation exact functors of
characteristic function $\rho \left( t\right) =t^{\theta }$. Then $\mathcal{F%
}_{\theta }\left( \overline{A}\right) =\left( \theta \left( 1-\theta \right) 
\frac{2}{2-\theta }\right) ^{\frac{2-\theta }{2}}.\left( A_{0},A_{1}\right)
_{\theta ,\frac{2}{2-\theta }}\footnote{%
For a Banach space $\left( X,\left\Vert {}\right\Vert \right) $ and $c>0$ we
note $c\cdot \left( X,\left\Vert {}\right\Vert \right) $ for the Banach
space $X$ with norm $c\left\Vert {}\right\Vert $.}$ and thus, for $\overline{%
A}=\left( L^{1}\left( 0,1\right) ,L^{2}\left( 0,1\right) \right) $ we get $%
\mathcal{F}_{\theta }\left( \overline{A}\right) =\left( \theta \left(
1-\theta \right) \frac{2}{2-\theta }\right) ^{\frac{2-\theta }{2}}.\left(
L^{1},L^{2}\right) _{\theta ,\frac{2}{2-\theta }}$. Thus $\mathcal{F}%
_{\theta }\left( \overline{A}\right) $ is $L^{\frac{2}{2-\theta }}=L^{p}$
with norm $\left( \theta \left( 1-\theta \right) \frac{2}{2-\theta }\right)
^{\frac{2-\theta }{2}}\phi _{\theta ,\frac{2}{2-\theta }}\left( K\left(
t,f,L^{1},L^{2}\right) \right) $ which is equivalent, for a fixed $p$, to
the usual norm of $L^{p}$, but having characteristic function $\rho \left(
t\right) =t^{\theta }$. We note $L^{p\blacktriangleleft }\left( 0,1\right)
:=\left( \theta \left( 1-\theta \right) \frac{2}{2-\theta }\right) ^{\frac{%
2-\theta }{2}}.\left( L^{1},L^{2}\right) _{\theta ,\frac{2}{2-\theta }}$ (or
simply $L^{p\blacktriangleleft }$), the space $L^{p}\left( 0,1\right) $ with
the norm 
\begin{equation*} 
\left\Vert f\right\Vert _{L^{p\blacktriangleleft }}:=\left( \theta \left(
1-\theta \right) \frac{2}{2-\theta }\right) ^{\frac{2-\theta }{2}}\phi
_{\theta ,\frac{2}{2-\theta }}\left( K\left( t,f,L^{1},L^{2}\right) \right) 
\end{equation*}%
. In this way we can use the following extrapolation result embedded in the
proof of Theorem 2.5 of \cite{JM} (see \cite{JM}, formulae in line 9 from
page 15, or \cite{AL} Theorem 2.2): If $\left\{ \mathcal{F}_{\theta
}\right\} _{\theta \in \Theta }$ is a complete family of exact interpolation
methods of characteristic function $\rho _{\theta }$ and $\overline{A}$ is a
mutually closed Banach pair, then $\left\Vert f\right\Vert _{\Delta \left(
\rho _{\theta }\left( t\right) \mathcal{F}_{\theta }\left( \overline{A}%
\right) \right) }\approx J\left( t,f,\overline{A}\right) $ with constants
independent of $t>0$. In particular, for $\mathcal{F}_{\theta }\left( 
\overline{A}\right) =L^{p\blacktriangleleft }=\left( \theta \left( 1-\theta
\right) \frac{2}{2-\theta }\right) ^{\frac{2-\theta }{2}}.\left(
L^{1},L^{2}\right) _{\theta ,\frac{2}{2-\theta }}$ we get 
\begin{equation*} 
\left\Vert f\right\Vert _{\Delta \left( t^{\theta }L^{p\blacktriangleleft
}\left( 0,1\right) \right) }\approx J\left( t,f,L^{1},L^{2}\right) 
\end{equation*}
and for $t=1$ this reads as $\left\Vert f\right\Vert _{\Delta \left(
L^{p\blacktriangleleft }\left( 0,1\right) \right) }\approx \max \left\{
\left\Vert f\right\Vert _{1},\left\Vert f\right\Vert _{2}\right\}
=\left\Vert f\right\Vert _{2}$ because $\left\Vert f\right\Vert _{1}\leq
\left\Vert f\right\Vert _{2}$ for all $f\in L^{1}\left( 0,1\right) $. With
all these facts at hand the equivalence between ii) and iii) in theorem \ref%
{N-B equivalence} is given from observing that $\left\Vert \chi _{\left(
0,1\right) }-f_{n}\right\Vert _{2}<\varepsilon $ for some $f_{n}\in \Phi $
and for any $\varepsilon >0$ if and only if $\sup\limits_{1<p<2}\left\Vert
\chi _{\left( 0,1\right) }-f_{n}\right\Vert _{L^{p\blacktriangleleft }\left(
0,1\right) }=\left\Vert \chi _{\left( 0,1\right) }-f_{n}\right\Vert _{\Delta
\left( L^{p\blacktriangleleft }\left( 0,1\right) \right) }<k\varepsilon $
for some $k>0$ independent from $f$, if and only if $\left\Vert \chi
_{\left( 0,1\right) }-f_{n}\right\Vert _{L^{p\blacktriangleleft }\left(
0,1\right) }<k\varepsilon $ for all $p\in \left( 1,2\right) $. In any case
let's remark that, by the equivalence of the norms for any fixed $p$, the
density of $\Phi $ in $L^{p}\left( 0,1\right) $ is equivalent to the density
of $\Phi $ in $L^{p\blacktriangleleft }\left( 0,1\right) $ for $1\leq p\leq
2 $.
\end{remark}

\begin{remark}
In addition, we also want to point out an equivalence for the $K-functional$
of the pair $\left( L^{1},L^{2}\right) $ which allows to give it in terms of
the well known $K-functional$ of the pair $\left( L^{1},L^{\infty }\right) $%
. We recall that if $K\left( t,f,\overline{X}\right) $ is the $K-functional$
for the couple $\overline{X}=\left( X_{0},X_{1}\right) $ then for $0<\theta
<1$, $1\leq q\leq \infty $ and $\overline{X}_{\theta ,q}=\left(
X_{0},X_{1}\right) _{\theta ,q}$ then Holmstedt's formulae for interpolation
spaces obtained by reiteration (cf. \cite{BL} Theorem 3.5.3) implies that 
\begin{equation} 
K\left( t,f,X_{0},\overline{X}_{\theta _{1},q}\right) \approx t\left(
\int\nolimits_{t^{1/\theta _{1}}}^{\infty }\left( s^{-\theta _{1}}K\left(
t,f,\overline{X}\right) \right) ^{q}\frac{ds}{s}\right) ^{\frac{1}{q}}
\label{Hoc} 
\end{equation}%
(see \cite{BL} corollary 3.6.2 i) pag 53, or \cite{BS} corollary V.2.3,
formula 2.19, pag 310). So, for the case $\overline{X}=\left(
L^{1},L^{\infty }\right) $ and $\theta _{1}=\frac{1}{2}$ and $q=2$, that is $%
\overline{X}_{\theta _{1},q}=\overline{X}_{\frac{1}{2},2}=L^{2}$ the
expression \ref{Hoc} reads like this: 
\begin{equation} 
K\left( t,f,L^{1},L^{2}\right) \approx t\left( \int\nolimits_{t^{2}}^{\infty
}\left( s^{-\frac{1}{2}}K\left( s,f,L^{1},L^{\infty }\right) \right) ^{2}%
\frac{ds}{s}\right) ^{\frac{1}{2}}  \label{L12a} 
\end{equation}%
and taking into account that $K\left( s,f,L^{1},L^{\infty }\right)
=\int_{0}^{s}f^{\ast }\left( u\right) du=sf^{\ast \ast }\left( s\right) $ we
can also write%
\begin{equation} 
K\left( t,f,L^{1},L^{2}\right) \approx t\left( \int\limits_{t^{2}}^{\infty
}\left( s^{\frac{1}{2}}f^{\ast \ast }\left( s\right) \right) ^{2}\frac{ds}{s}%
\right) ^{\frac{1}{2}}  \label{L12b} 
\end{equation}
\end{remark}

,that is 
\begin{equation} 
\frac{K\left( t,f,L^{1},L^{2}\right) }{t}\approx \left(
\int\limits_{t^{2}}^{\infty }\left( f^{\ast \ast }\left( s\right) \right)
^{2}ds\right) ^{\frac{1}{2}}  \label{L12c} 
\end{equation}%
.

\begin{remark}
Returning to the statement of the previous theorem \ref{N-B equivalence} the
lack of zeros of $\zeta $ in $\left( \frac{1}{2},1\right) $ is equivalent to
the lack of zeros in $\left( \frac{1}{p},1\right) $ for all $p\in \left(
1,2\right) $ that's why iii) implies ii); the converse, on the other hand,
is obvious.
\end{remark}

So, one way to prove the Riemann Hypothesis, if it were true, could be by
proving the density of $\Phi $ in $L^{p}\left( 0,1\right) $ for every $p\in
\left( 1,2\right) $.

Now we come to the foremost matter, that involve to replace the
interpolation spaces $\left( L^{1},L^{2}\right) _{\theta ,\frac{1}{1-\frac{%
\theta }{2}}}$ for $\theta \in \left( 0,1\right) $ which we will give us $%
L^{p}$ with $p\in \left( 1,2\right) $ by the larger spaces $M\left( \theta
\right) \left( L^{1},L^{2}\right) _{\theta ,\infty ,K}$ improving the above
mentioned requirement about $L^{p}$, from Nyman-Beurling theorem, for
proving, if true, the Riemann Hypothesis by means of weaker conditions:

We begin by recalling one more theorem of Jawerth-Milman Extrapolation
Theory, we copy its statement from \cite{M} (Theorem 21, pag 44 and (4.1) in
the preceding remarks)

\begin{theorem}
\label{M_delta}Let $\overline{A}=\left( A_{0},A_{1}\right) $ a mutually
closed Banach couple, and let $M\left( \theta \right) $ a tempered function
on $\Theta =\left( 0,1\right) $. For any $\left\{ \mathcal{F}_{\theta
}\right\} _{\theta \in \Theta }$ a family of interpolation exact functors of
order $\theta $, i.e. with characteristic function $t^{\theta }$, we have
that: 
\begin{equation*} 
\Delta _{0<\theta <1}\left( M\left( \theta \right) \mathcal{F}_{\theta
}\left( \overline{A}\right) \right) =\Delta _{0<\theta <1}\left( M\left(
\theta \right) A_{\theta ,\infty ,K}\right) 
\end{equation*}%
with equivalent norms.

Moreover if $A_{1}\subset A_{0}$, which means that $\overline{A}$ is an
ordered pair, we have 
\begin{equation} 
\Delta _{\theta _{0}<\theta <1}\left( M\left( \theta \right) \mathcal{F}%
_{\theta }\left( \overline{A}\right) \right) =\Delta _{0<\theta <1}\left(
M\left( \theta \right) \mathcal{F}_{\theta }\left( \overline{A}\right)
\right)  \label{01eqtheta01} 
\end{equation}%
for any fixed $\theta _{0}\in \left( 0,1\right) $ and then 
\begin{equation*} 
\Delta _{\theta _{0}<\theta <1}\left( M\left( \theta \right) \mathcal{F}%
_{\theta }\left( \overline{A}\right) \right) =\Delta _{\theta _{0}<\theta
<1}\left( M\left( \theta \right) A_{\theta ,\infty ,K}\right) 
\end{equation*}
\end{theorem}

Now it's time to use the Jawerth-Milman theorem \ref{M_delta} for the $\Delta
-functor$ of extrapolation: We consider $p=\frac{1}{1-\frac{\theta }{2}}$
where $p\in \left( 1,2\right) $ while $\theta \in \left( 0,1\right) $ and,
as we saw in a previous paragraph, $L^{p}=L^{\frac{1}{1-\frac{\theta }{2}}%
}=\left( L^{1},L^{2}\right) _{\theta ,\frac{1}{1-\frac{\theta }{2}},K}$ for $%
K=K\left( t,f,L^{1},L^{2}\right) $, then if we consider $L^{p%
\blacktriangleleft }=L^{\frac{1}{1-\frac{\theta }{2}}\blacktriangleleft }$
by renorming $L^{p}$ with the norm $\left\Vert \cdot \right\Vert
_{L^{p\blacktriangleleft }}:=\left( p\theta \left( 1-\theta \right) \right)
^{\frac{1}{p}}\left\Vert \cdot \right\Vert _{L^{p}}$, with $p=\frac{1}{1-%
\frac{\theta }{2}}$ and understanding $\left( p\theta \left( 1-\theta
\right) \right) ^{\frac{1}{p}}$ as $1$ if $p=\infty $, then the
interpolation functor has characteristic function $\rho _{\theta }\left(
t\right) =t^{\theta }$, so $\rho ^{\ast }\left( t\right)
=\inf\limits_{\theta \in \left( 0,1\right) }\rho _{\theta }\left( t\right)
=\min \{1,t\}$. Now for a tempered $M\left( \theta \right) $, by the theorem %
\ref{M_delta} we have that 
\begin{eqnarray*} 
\Delta _{0<\theta <1}\left( M\left( \theta \right) L^{\frac{1}{1-\frac{%
\theta }{2}}\blacktriangleleft }\left( 0,1\right) \right) &=&\Delta
_{0<\theta <1}\left( M\left( \theta \right) \mathcal{F}_{\theta }\left( 
\overline{A}\right) \right) \\
&=&\Delta _{0<\theta <1}\left( M\left( \theta \right) \overline{A}_{\theta
,\infty ,K}\right) 
\end{eqnarray*}%
. Therefore for every $f\in L^{1}:$

\begin{eqnarray} 
\left\Vert f\right\Vert _{\Delta _{0<\theta <1}\left( M\left( \theta \right)
L^{\frac{1}{1-\frac{\theta }{2}}\blacktriangleleft }\left( 0,1\right)
\right) } &\approx &\left\Vert f\right\Vert _{\Delta _{0<\theta <1}\left(
M\left( \theta \right) \overline{A}_{\theta ,\infty ,K}\right) }
\label{delta_L1L2} \\ 
&=&\sup\limits_{0<\theta <1}\left\Vert f\right\Vert _{M\left( \theta \right) 
\overline{A}_{\theta ,\infty ,K}}  \notag \\ 
&=&\sup\limits_{0<\theta <1}(\sup\limits_{t>0}M\left( \theta \right) \frac{%
K\left( t,f,L^{1},L^{2}\right) }{t^{\theta }})  \notag 
\end{eqnarray}

Now using that $K\left( t,f,L^{1},L^{2}\right) \approx t\left(
\int\limits_{t^{2}}^{\infty }\left( s^{\frac{1}{2}}f^{\ast \ast }\left(
s\right) \right) ^{2}\frac{ds}{s}\right) ^{\frac{1}{2}}$ (cf. \ref{L12b}
above) we have:%
\begin{equation} 
\left\Vert f\right\Vert _{\Delta _{_{0<\theta <1}}\left( M\left( \theta
\right) L_{\left( 0,1\right) }^{\frac{1}{1-\frac{\theta }{2}}%
\blacktriangleleft }\right) }\approx \sup\limits_{0<\theta <1}\left(
\sup\limits_{t>0}M\left( \theta \right) t^{1-\theta }\left(
\int\limits_{t^{2}}^{\infty }f^{\ast \ast }\left( s\right) ^{2}ds\right) ^{%
\frac{1}{2}}\right)  \label{delta_theta_inf} 
\end{equation}%
. Let's remark that the constants of the above equivalence \ref{delta_theta_inf}
don't depend on $\theta $ because we only use the equivalence \ref{L12c} for the 
$K-functional$ for the pair $\left( L^{1},L^{2}\right) $.

Now, for each $\theta \in \left( 0,1\right) $ is $p=\frac{1}{1-\frac{\theta 
}{2}}\in \left( 1,2\right) $ and $t^{1-\theta }=t^{1-\frac{2}{p^{\prime }}}$
and putting $\omega \left( p\right) =M\left( \theta \right) $ for any $f$
and writing $L^{p\blacktriangleleft }$ instead of $L^{p\blacktriangleleft
}\left( 0,1\right) $ we have that (\ref{delta_theta_inf}) reads

\begin{eqnarray} 
\left\Vert f\right\Vert _{\Delta _{1<p<2}\left( \omega \left( p\right)
L^{p\blacktriangleleft }\right) } &=&\left\Vert f\right\Vert _{\Delta
_{0<\theta <1}\left( M\left( \theta \right) L^{\frac{1}{1-\frac{\theta }{2}}%
\blacktriangleleft }\right) }  \label{delta_LpInf} \\ 
&\approx &\sup\limits_{1<p<2}\left( \sup\limits_{t>0}\omega \left( p\right)
\left( \frac{2}{p^{\prime }}\right) ^{\frac{2}{p^{\prime }}-1}t^{1-\frac{2}{%
p^{\prime }}}\left( \int\limits_{t^{2}}^{\infty }f^{\ast \ast }\left(
s\right) ^{2}ds\right) ^{\frac{1}{2}}\right)  \notag 
\end{eqnarray}

Bear in mind that it is easy to prove that $\left\Vert \cdot \right\Vert
_{\Delta _{0<\theta <1}\left( M\left( \theta \right) L^{\frac{1}{1-\frac{%
\theta }{2}}\blacktriangleleft }\left( 0,1\right) \right) }=\left\Vert \cdot
\right\Vert _{\Delta _{\theta _{0}<\theta <1}\left( M\left( \theta \right)
L^{\frac{1}{1-\frac{\theta }{2}}\blacktriangleleft }\left( 0,1\right)
\right) }$ for any fixed $\theta _{0}\in \left( 0,1\right) $ (cf. \cite{M}
(4.1) pag. 44) we can consider just the supremum in $p_{0}<p<2$ and observing
that $\left( \frac{2}{p^{\prime }}\right) ^{\frac{2}{p^{\prime }}%
-1}\rightarrow 1$ when $p$ goes to $2$ we can drop the factor $\left( \frac{2%
}{p^{\prime }}\right) ^{\frac{2}{p^{\prime }}-1}$ and we have that:%
\begin{eqnarray} 
&&\left\Vert f\right\Vert _{\Delta _{p_{0}<p<2}\left( \omega \left( p\right)
L^{p\blacktriangleleft }\left( 0,1\right) \right) }  \label{ineq_AthetaInf}
\\
&\approx &\left\Vert f\right\Vert _{\Delta _{1<p<2}\left( \omega \left(
p\right) L^{p\blacktriangleleft }\left( 0,1\right) \right) }  \notag \\
&\approx &\sup\limits_{1<p<2}\left( \sup\limits_{t>0}\omega \left( p\right)
t^{1-\frac{2}{p^{\prime }}}\left( \int\limits_{t^{2}}^{\infty }f^{\ast \ast
}\left( s\right) ^{2}ds\right) ^{\frac{1}{2}}\right)  \notag 
\end{eqnarray}

Taking into account that, for the probability space $\left( 0,1\right) $,
the pair

$\overline{A}=\left( L^{1},L^{2}\right) $ is ordered we can use the
equivalence between the usual interpolation spaces $\overline{A}_{\theta
,q,K}$ spaces and the modified ones $\left\langle \overline{A}_{\theta
,q,K}\right\rangle =\{f\in A_{0}:\left\Vert f\right\Vert _{\left\langle 
\overline{A}_{\theta ,q,K}\right\rangle }<\infty \}\ $with equivalence of
norms (cf. \cite{ALM} proposition 1 section 3.2): 
\begin{equation} 
\left\Vert f\right\Vert _{\left\langle \overline{A}_{\theta
,q,K}\right\rangle }\leq \left\Vert f\right\Vert _{\overline{A}_{\theta
,q,K}}\leq \left( 1+\left( 1-\theta \right) ^{\frac{1}{q}}\theta ^{-\frac{1}{%
q}}\right) \left\Vert f\right\Vert _{\left\langle \overline{A}_{\theta
,q,K}\right\rangle }  \label{Modified spaces equivalence} 
\end{equation}%
, where $\left\Vert f\right\Vert _{\left\langle \overline{A}_{\theta
,q,K}\right\rangle }=\phi _{\theta ,q}\left( \chi _{\left( 0,1\right)
}K\left( s,f,\overline{A}\right) \right) $, that is 
\begin{equation*} 
\left\Vert f\right\Vert _{\left\langle \overline{A}_{\theta
,q,K}\right\rangle }=\left\{ 
\begin{array}{c}
\left( \int_{0}^{1}\left( s^{-\theta }K\left( s,f,\overline{A}\right)
\right) ^{q}\frac{ds}{s}\right) ^{\frac{1}{q}}\text{ if }q<\infty \\ 
\sup\limits_{0<s<1}s^{-\theta }K\left( s,f,\overline{A}\right) \text{ if }%
q=\infty%
\end{array}%
\right. 
\end{equation*}%
. Also, in the same way of the normalizing explained in a remark above, we
can consider $\left\langle \overline{A}_{\theta ,q,K}\right\rangle
^{\blacktriangleleft }:=\left( \theta \left( 1-\theta \right) \frac{2}{%
2-\theta }\right) ^{\frac{2-\theta }{2}}.\left\langle \overline{A}_{\theta
,q,K}\right\rangle $ to obtain characteristic functions that are exactly $%
\rho _{\theta }\left( t\right) =t^{\theta }$.

Although in \ref{Modified spaces equivalence} the factor blows up when $%
\theta $ goes to zero the endpoint what matters appears for $\theta
\rightarrow 1$ where $\left( 1+\left( 1-\theta \right) ^{\frac{1}{q}}\theta
^{-\frac{1}{q}}\right) $ is bounded, so by fixing some $\theta _{0}\in
\left( 0,1\right) $ we have that $\left\Vert f\right\Vert _{\left\langle 
\overline{A}_{\theta ,q,K}\right\rangle }\approx
\sup\limits_{0<s<1}s^{-\theta }K\left( s,f,\overline{A}\right) $ for $\theta
\in \left( \theta _{0},1\right) $, and recalling that $\Delta _{\theta
_{0}<\theta <1}\left( M\left( \theta \right) \mathcal{F}_{\theta }\left( 
\overline{A}\right) \right) =\Delta _{0<\theta <1}\left( M\left( \theta
\right) \mathcal{F}_{\theta }\left( \overline{A}\right) \right) $ we arrive
at 
\begin{equation*} 
\left\Vert f\right\Vert _{\Delta _{_{0<\theta <1}}\left( M\left( \theta
\right) L^{\frac{1}{1-\frac{\theta }{2}}\blacktriangleleft }\left(
0,1\right) \right) }\approx \sup\limits_{0<\theta <1}\left(
\sup\limits_{0<t<1}M\left( \theta \right) t^{1-\theta }\left(
\int\limits_{t^{2}}^{\infty }f^{\ast \ast }\left( s\right) ^{2}ds\right) ^{%
\frac{1}{2}}\right) 
\end{equation*}%
and also 
\begin{equation*} 
\left\Vert f\right\Vert _{\Delta _{1<p<2}\left( \omega \left( p\right)
L^{p\blacktriangleleft }\left( 0,1\right) \right) }\approx
\sup\limits_{1<p<2}\left( \sup\limits_{0<t<1}\omega \left( p\right) t^{1-%
\frac{2}{p^{\prime }}}\left( \int\limits_{t^{2}}^{\infty }f^{\ast \ast
}\left( s\right) ^{2}ds\right) ^{\frac{1}{2}}\right) 
\end{equation*}%
, and thus we replace the supremums of $t$ for $t>0$ by supremums of $t$ for $%
0<t<1$.

Now we are ready to see the main result and its different formulations in
terms of the $K-functional$, as a direct application of the Jawerth-Milman
theorem \ref{M_delta} and the Nyman-Beurling theorem \ref{N-B equivalence}.
Our result improves the criterion given by theorem \ref{N-B equivalence}. Let's recall the statement that we advanced above:

\begin{criterion}
\label{N-B-with-extrapol}The following statements are equivalent:

a) Riemann Hypothesis is true

b) For $\overline{A}=\left( L^{1}\left( 0,1\right) ,L^{2}\left( 0,1\right)
\right) $ there is a positive function $M\left( \theta \right) $ tempered at 
$1$ with $M\left( \theta \right) =O(1)\footnote{%
We get a real improvement respect to the Nyman-Beurling Criterion when $%
M\left( \theta \right) =o(1)$ for $\theta \rightarrow 1$}$ for $\theta
\rightarrow 1$ such that%
\begin{equation*} 
\inf\limits_{f\in \Phi }\left\Vert \chi _{\left( 0,1\right) }-f\right\Vert
_{\Delta _{0<\theta <1}\left( M\left( \theta \right) \overline{A}_{\theta
,\infty ,K\left( \cdot ,\cdot ,L^{1},L^{2}\right) }\right) }=0 
\end{equation*}

That is, 
\begin{equation*} 
\sup\limits_{0<\theta <1}\left( \inf\limits_{f\in \Phi }\left[
\sup\limits_{0<t<1}M\left( \theta \right) t^{1-\theta }\left(
\int\limits_{t^{2}}^{\infty }\left( \chi _{\left( 0,1\right) }-f\right)
^{\ast \ast }\left( s\right) ^{2}ds\right) ^{\frac{1}{2}}\right] \right) =0 
\end{equation*}
\end{criterion}

\begin{proof}
The implication $a)\Longrightarrow b)$ follows easily from the
Nyman-Beurling theorem (\ref{N-B equivalence}): Since RH implies that $\Phi $
is dense in $L^{2\blacktriangleleft }$ by \ref{N-B equivalence}, and being 
\begin{equation*} 
\left\Vert \cdot \right\Vert _{\overline{A}_{\theta ,\infty ,K\left( \cdot
,\cdot ,L^{1},L^{2}\right) }}\leq \left\Vert \cdot \right\Vert _{\overline{A}%
_{\theta ,\frac{2}{2-\theta },K\left( \cdot ,\cdot ,L^{1},L^{2}\right)
}}=\left\Vert \cdot \right\Vert _{L^{\frac{2}{2-\theta }}\left( 0,1\right)
}\leq \left\Vert \cdot \right\Vert _{L^{2}\left( 0,1\right) }\approx
\left\Vert \cdot \right\Vert _{L^{2\blacktriangleleft }\left( 0,1\right) } 
\end{equation*}
for $\theta \in \left( 0,1\right) $, then given any $\varepsilon >0$ there
is $f\in \Phi $ such that $\left\Vert \chi _{\left( 0,1\right)
}-f\right\Vert _{L^{2\blacktriangleleft }\left( 0,1\right) }<\varepsilon $.
Therefore, because $M\left( \theta \right) =O(1)$ for $\theta \rightarrow 1$%
, we can get $\theta _{0}<1$ such that for some $k>0$ and for any $\theta
\in \left( \theta _{0},1\right) $ 
\begin{eqnarray*} 
\left\Vert \chi _{\left( 0,1\right) }-f\right\Vert _{M\left( \theta \right) 
\overline{A}_{\theta ,\infty ,K\left( \cdot ,\cdot ,L^{1},L^{2}\right) }}
&=&M\left( \theta \right) \left\Vert \chi _{\left( 0,1\right) }-f\right\Vert
_{\overline{A}_{\theta ,\infty ,K\left( \cdot ,\cdot ,L^{1},L^{2}\right) }}
\\
&\leq &k\left\Vert \chi _{\left( 0,1\right) }-f\right\Vert _{\overline{A}%
_{\theta ,\frac{2}{2-\theta },K\left( \cdot ,\cdot ,L^{1},L^{2}\right) }} \\
&=&k\left\Vert \chi _{\left( 0,1\right) }-f\right\Vert _{L^{\frac{2}{%
2-\theta }}\left( 0,1\right) } \\
&\leq &k\left\Vert \chi _{\left( 0,1\right) }-f\right\Vert _{L^{2}\left(
0,1\right) }<\varepsilon 
\end{eqnarray*}

Thus renaming $\varepsilon =\frac{\varepsilon }{k}:$ 
\begin{equation} 
\left\Vert \chi _{\left( 0,1\right) }-f\right\Vert _{\Delta _{\theta
_{0}<\theta <1}\left( M\left( \theta \right) \overline{A}_{\theta ,\infty
,K\left( \cdot ,L^{1},L^{2}\right) }\right) }\leq \varepsilon 
\label{inf_theta0} 
\end{equation}%
\footnote{%
Let's remember that $\overline{X}_{\theta ,\infty ,K}^{\blacktriangleleft }=%
\overline{X}_{\theta ,\infty ,K}$, so we write $\overline{A}_{\theta ,\infty
,K\left( \cdot ,L^{1},L^{2}\right) }$ instead of $\overline{A}_{\theta
,\infty ,K\left( \cdot ,L^{1},L^{2}\right) }^{\blacktriangleleft }.$}.

But from the equivalence \ref{01eqtheta01} in the statement of the
Jawerth-Milman theorem \ref{M_delta} we have that 
\begin{equation*} 
\Delta _{\theta _{0}<\theta <1}\left( M\left( \theta \right) \overline{A}%
_{\theta ,\infty ,K\left( \cdot ,\cdot ,L^{1},L^{2}\right) }\right) =\Delta
_{0<\theta <1}\left( M\left( \theta \right) \overline{A}_{\theta ,\infty
,K\left( \cdot ,\cdot ,L^{1},L^{2}\right) }\right) 
\end{equation*}%
with equivalence of the norms. So from inequality \ref{inf_theta0}, for any $%
\varepsilon >0$ we can find $f\in \Phi $ such that: 
\begin{equation*} 
\left\Vert \chi _{\left( 0,1\right) }-f\right\Vert _{\Delta _{0<\theta
<1}\left( M\left( \theta \right) \overline{A}_{\theta ,\infty ,K\left( \cdot
,\cdot ,L^{1},L^{2}\right) }\right) }<\varepsilon 
\end{equation*}%
and then $\inf\limits_{f\in \Phi }\left\Vert \chi _{\left( 0,1\right)
}-f\right\Vert _{\Delta _{0<\theta <1}\left( M\left( \theta \right) 
\overline{A}_{\theta ,\infty ,K\left( \cdot ,\cdot ,L^{1},L^{2}\right)
}\right) }=0$, so we have arrived at $b)$.

Now let's assume $b)$: $\inf\limits_{f\in \Phi }\left\Vert \chi _{\left( 0,1\right)
}-f\right\Vert _{\Delta _{0<\theta <1}\left( M\left( \theta \right) 
\overline{A}_{\theta ,\infty ,K\left( \cdot ,\cdot ,L^{1},L^{2}\right)
}\right) }=0$.

Taking in account that $\overline{A}_{\theta ,\frac{2}{2-\theta },K\left(
\cdot ,\cdot ,L^{1}\left( 0,1\right) ,L^{2}\left( 0,1\right) \right) }=L^{%
\frac{2}{2-\theta }}\left( 0,1\right) $ whose norm is equivalent to $%
\left\Vert \cdot \right\Vert _{L^{\frac{2}{2-\theta }\blacktriangleleft
}\left( 0,1\right) }$ for any fixed $\theta $ and using the main statement
of Jawerth-Milman theorem \ref{M_delta} we have that 
\begin{equation*} 
\Delta _{0<\theta <1}\left( M\left( \theta \right) \overline{A}_{\theta
,\infty ,K\left( \cdot ,\cdot ,L^{1},L^{2}\right) }\right) =\Delta
_{0<\theta <1}\left( M\left( \theta \right) L^{\frac{2}{2-\theta }%
\blacktriangleleft }\left( 0,1\right) \right) 
\end{equation*}%
with equivalence of norms. Thus there is $f\in \Delta _{0<\theta <1}\left(
M\left( \theta \right) L^{\frac{2}{2-\theta }\blacktriangleleft }\left(
0,1\right) \right) $ with $\left\Vert \chi _{\left( 0,1\right)
}-f\right\Vert _{\Delta _{0<\theta <1}\left( M\left( \theta \right) L^{\frac{%
2}{2-\theta }\blacktriangleleft }\left( 0,1\right) \right) }<\varepsilon $,
for any $\varepsilon >0$. Now, again for any $\delta >0$ we can choose $%
\delta _{1}>0$ with $0<\delta _{1}\leq \delta $ small enough to have that if
we take $p=2-\delta _{1}=\frac{2}{2-\theta }$ then $\theta =2-\frac{2}{p}\in
\left( \theta _{0},1\right) $, where again $\theta _{0}$ is close enough to $%
1$ to ensure that $M\left( \theta \right) \leq k$ for some $k>0$ and for all 
$\theta \in \left( \theta _{0},1\right) $. Hence we obtain that $k\left\Vert
\chi _{\left( 0,1\right) }-f\right\Vert _{L^{p}\left( 0,1\right) }\lesssim
k\left\Vert \chi _{\left( 0,1\right) }-f\right\Vert _{L^{p\blacktriangleleft
}\left( 0,1\right) }\leq M\left( \theta \right) \left\Vert \chi _{\left(
0,1\right) }-f\right\Vert _{L^{\frac{2}{2-\theta }\blacktriangleleft }\left(
0,1\right) }<\varepsilon $. Therefore we have that $\chi _{\left( 0,1\right)
}\in \overline{\Phi }^{L^{p}}$ and then Nyman-Beurling theorem \ref%
{Nyman-Beurling} guarantees that $\zeta \left( s\right) $ is zero-free
in the half-plane $\sigma >\frac{1}{p}$, where $\sigma =\operatorname{Re}\left(
s\right) $, but we can take $p$ so close to $2$ as we like and then $\zeta
\left( s\right) $ is zero-free in the half-plane $\sigma >\frac{1}{2}$
and the Riemann Hypothesis would hold. So $b)$ implies $a)$.
\end{proof}

\begin{remark}
The above result can be written in terms of $p\in \left( 1,2\right) $ such
that $\theta =\frac{2}{p^{\prime }}=2-\frac{2}{p}$ and requiring for a
function $\omega \left( p\right) $ such that 
\begin{equation*} 
\omega \left( p\right) \left( L^{1},L^{2}\right) _{\frac{2}{p^{\prime }}%
,\infty ,K\left( \cdot ,\cdot ,L^{1},L^{2}\right) }=M\left( \theta \right)
\left( L^{1},L^{2}\right) _{\theta ,\infty ,K\left( \cdot ,\cdot
,L^{1},L^{2}\right) } 
\end{equation*}
with equivalent norms then the equality of statement b) reads like this:

b')%
\begin{equation*} 
\inf\limits_{f\in \Phi }\left\Vert \chi _{\left( 0,1\right) }-f\right\Vert
_{\Delta _{1<p<2}\left( \omega \left( p\right) \left( L^{1},L^{2}\right) _{%
\frac{2}{p^{\prime }},\infty ,K\left( \cdot ,\cdot ,L^{1},L^{2}\right)
}\right) }=0 
\end{equation*}%
Therefore the conditions of b) and b') respectively mean that:%
\begin{equation*} 
\sup\limits_{0<\theta <1}\left( \inf\limits_{f\in \Phi }\left[
\sup\limits_{0<t<1}M\left( \theta \right) t^{1-\theta }\left(
\int\limits_{t^{2}}^{\infty }\left( \chi _{\left( 0,1\right) }-f\right)
^{\ast \ast }\left( s\right) ^{2}ds\right) ^{\frac{1}{2}}\right] \right) =0 
\end{equation*}%
or 
\begin{equation*} 
\sup\limits_{1<p<2}\left( \inf\limits_{f\in \Phi }\left[ \sup\limits_{0<t<1}%
\omega \left( p\right) t^{1-\frac{2}{p^{\prime }}}\left(
\int\limits_{t^{2}}^{\infty }\left( \chi _{\left( 0,1\right) }-f\right)
^{\ast \ast }\left( s\right) ^{2}ds\right) ^{\frac{1}{2}}\right] \right) =0 
\end{equation*}
\end{remark}

\begin{problem}
Taking in account that Grand Lebesgue spaces $L^{p)}$ can be seen as $\Delta
-$extrapolation spaces and that the quantities in b')\ look a bit analogous
to Fiorenza-Karadzhov's equivalence to Grand Lebesgue $L^{p)}$ norms (cf. 
\cite{FK}) it seems possible to obtain an analogous criterion for RH in
terms of the density of $\Phi $ in $L^{p)}$ spaces.
\end{problem}

Writing $K=K\left( t,f,\overline{A}\right) $ for $\overline{A}=\left(
L^{1},L^{2}\right) $ in terms of $K\left( t,f,L^{1},L^{\infty }\right) $ by
using reiteration we can write the criterion \ref{N-B-with-extrapol} in a
terse form:

\begin{theorem}
\label{strong Nyman-Beurling-Milman equivalence}The following statements are
equivalent:

i) \textit{Riemann zeta-function }$\zeta $\textit{\ is zero-free for }$%
\frac{1}{2}<\operatorname{Re}\left( s\right) <1$, or in other words the Riemann
Hypothesis is true.

ii) $\chi _{\left( 0,1\right) }\in \overline{\Phi }^{\Delta _{0<\theta
<1}\left( M\left( \theta \right) \overline{A}_{\theta ,\infty ,K\left( \cdot
,\cdot ,L^{1},L^{2}\right) }\right) }$.
\end{theorem}

\begin{remark}
Taking $p=\frac{2}{2-\theta }$ and using that $\left( L^{1},L^{2}\right) _{%
\frac{2}{p^{\prime }},p,K}\overset{1}{\hookrightarrow }\left(
L^{1},L^{2}\right) _{\frac{2}{p^{\prime }},\infty ,K}$ we have that $%
\left\Vert \chi _{\left( 0,1\right) }-f\right\Vert _{\omega \left( p\right)
\left( L^{1},L^{2}\right) _{\frac{2}{p^{\prime }},\infty ,K}}\leq \left\Vert
\chi _{\left( 0,1\right) }-f\right\Vert _{\omega \left( p\right) L^{p}}$. We
have shown that to prove the RH it is enough to find, for any $\varepsilon
>0 $ and any $p\in \left( 1,2\right) $, some $f\in \Phi $ such that $%
\sup\limits_{1<\theta <2}\left\Vert \chi _{\left( 0,1\right) }-f\right\Vert
_{\omega \left( p\right) \left( L^{1},L^{2}\right) _{\frac{2}{p^{\prime }}%
,\infty ,K}}<\varepsilon $, or in other words for some fixed $C>0$: 
\begin{equation*} 
\sup\limits_{1<p<2}\left( \sup\limits_{0<t<1}\omega \left( p\right) t^{1-%
\frac{2}{p^{\prime }}}\left( \int\limits_{t^{2}}^{\infty }\left( \chi
_{\left( 0,1\right) }-f\right) ^{\ast \ast }\left( s\right) ^{2}ds\right) ^{%
\frac{1}{2}}\right) <C\varepsilon 
\end{equation*}%
.

For a function $M\left( \theta \right) $ tempered at $1$ such that $M\left(
\theta \right) \rightarrow 0$ for $\theta \rightarrow 1$, or respectively $%
\omega \left( p\right) \rightarrow 0$ for $p\rightarrow 2$ this requirement
is in general weaker to the one of the Nyman-Beurling original criterion,
which can be read as $\sup\limits_{1<p<2}\left( \left\Vert \chi _{\left(
0,1\right) }-f\right\Vert _{L^{p}}\right) <\varepsilon $, for any $%
\varepsilon >0$ and for any $p\in \left( 1,2\right) $\footnote{%
It's clear that, for any $p$, $\inf\limits_{1<p<2}\omega \left( p\right)
\left\Vert \chi _{\left( 0,1\right) }-f\right\Vert
_{L^{p}}=\inf\limits_{1<p<2}\left\Vert \chi _{\left( 0,1\right)
}-f\right\Vert _{\omega \left( p\right) L^{p}}=0$ if and only if $%
\inf\limits_{1<p<2}\left\Vert \chi _{\left( 0,1\right) }-f\right\Vert
_{L^{p}}=0$}. This occurs because $\omega \left( p\right) \left\Vert \chi
_{\left( 0,1\right) }-f\right\Vert _{\left( L^{1},L^{2}\right) _{\frac{2}{%
p^{\prime }},\infty ,K}}\leq \left\Vert \chi _{\left( 0,1\right)
}-f\right\Vert _{L^{p}}$ and so, clearly, $\inf\limits_{f\in \Phi
}\left\Vert \chi _{\left( 0,1\right) }-f\right\Vert _{L^{p}}=0$ implies $%
\inf\limits_{f\in \Phi }\left\Vert \chi _{\left( 0,1\right) }-f\right\Vert
_{\omega \left( p\right) \left( L^{1},L^{2}\right) _{\frac{2}{p^{\prime }}%
,\infty ,K}}=0$, but the converse is not always true, so for to prove the RH
we could find a family $\left\{ f_{\varepsilon }\right\} _{\varepsilon
>0}\subset $ $\Phi $ such that $\sup\limits_{1<p<2}\left( \left\Vert \chi
_{\left( 0,1\right) }-f_{\varepsilon }\right\Vert _{\omega \left( p\right)
\left( L^{1},L^{2}\right) _{\frac{2}{p^{\prime }},\infty ,K}}\right)
<\varepsilon $ even though it could happen that $\inf\limits_{\varepsilon
>0}\left\Vert \chi _{\left( 0,1\right) }-f_{\varepsilon }\right\Vert
_{L^{2}}\geq c>0$ $\footnote{%
Nevertheless by proving ii) we would've proved RH, giving an indirect proof
of $\inf\limits_{f\in \Phi }\left\Vert \chi _{\left( 0,1\right)
}-f\right\Vert _{L^{2}}=0$ by means of Nyman-Beurling theorem even if we
wouldn't have exhibit a family $\left\{ g_{\varepsilon }\right\}
_{\varepsilon >0}\subset \Phi $ such that $\left\Vert \chi _{\left(
0,1\right) }-g_{\varepsilon }\right\Vert _{L^{2}}<\varepsilon $}$. So, in
this sense our criterion is sharper to the one from Nyman-Beurling when you
try to prove the (hypothetical) trueness of Riemann Hypothesis.
\end{remark}

\begin{remark}
We wrote the above results in terms of the density of the family of
functions $\Phi =\left\{ \sum\limits_{k=1}^{n}c_{n_{k}}\rho \left( \frac{1}{%
a_{n_{k}}x}\right) :\sum\limits_{k=1}^{n}c_{n_{k}}\rho \left( \frac{1}{%
a_{n_{k}}}\right) =0\right\} _{a_{n_{k}}\geq 1}$ but we could have used
instead $\Phi _{0}=\left\{ \sum\limits_{k=1}^{n}c_{n_{k}}\rho \left( \frac{1%
}{n_{k}x}\right) :\sum\limits_{k=1}^{n}c_{n_{k}}\rho \left( \frac{1}{n_{k}}%
\right) =0\right\} _{n_{k}\in \mathbb{N}}$ taking into account that, as we
have already mentioned, L. B\'{a}ez-Duarte proved that $\overline{\Phi _{0}}%
^{L^{2}}=\overline{\Phi }^{L^{2}}$.
\end{remark}

\begin{remark}
A remarkable point is that actually we have a large family of criteria,
depending on which tempered function $M\left( \theta \right) $ is picked.
Some of them sharper than those in corollary \ref{NBC}. Let's also observe that for the constant function $M\left( \theta \right) \equiv
1$, which is trivially tempered, we recover the original Nyman-Beurling
criterion for the Riemann Hypothesis, let's see the instructive details:
\end{remark}

If $M\left( \theta \right) =1$ it is 
\begin{equation*} 
\Delta _{1<p<2}\left( \omega \left( p\right) L^{p\blacktriangleleft }\right)
=\Delta _{0<\theta <1}\left( L^{\frac{1}{1-\frac{\theta }{2}}%
\blacktriangleleft }\left( 0,1\right) \right) =\Delta _{0<\theta <1}\left( 
\overline{A}_{\theta ,\infty ,K}\right) =A_{\rho ^{\ast },\infty ,K} 
\end{equation*}%
where $\rho ^{\ast }\left( t\right) =\inf\limits_{\theta \in \left(
0,1\right) }\rho _{\theta }\left( t\right) =\inf\limits_{\theta \in \left(
0,1\right) }t^{\theta }=\min \{1,t\}$ being $\rho _{\theta }\left( t\right)
=t^{\theta }$ the characteristic function of the functors $\overline{A}%
_{\theta ,q,K}$. The last equality: $\Delta _{0<\theta <1}\left( \overline{A}%
_{\theta ,\infty ,K}\right) =A_{\rho ^{\ast },\infty ,K}$, due to Jawerth
and Milman, follows easily from a "Fubini argument" (cf. \cite{AL} theorem
2.11). Then using this and \ref{L12b} we have

\begin{eqnarray*} 
\left\Vert f\right\Vert _{A_{\rho ^{\ast },\infty ,K}} &=&\sup\limits_{t>0}%
\frac{K\left( t,f,L^{1},L^{2}\right) }{\min \{1,t\}}=\sup\limits_{t>0}\left(
\max \{1,\frac{1}{t}\}K\left( t,f,L^{1},L^{2}\right) \right) \\
&\approx &\sup\limits_{t>0}\left( \max \{1,\frac{1}{t}\}t\left(
\int\nolimits_{t^{2}}^{\infty }\left( s^{\frac{1}{2}}f^{\ast \ast }\left(
s\right) \right) ^{2}\frac{ds}{s}\right) ^{\frac{1}{2}}\right) \\
&=&\sup\limits_{t>0}\left( \max \{t,1\}\left( \int\nolimits_{t^{2}}^{\infty
}\left( s^{\frac{1}{2}}f^{\ast \ast }\left( s\right) \right) ^{2}\frac{ds}{s}%
\right) ^{\frac{1}{2}}\right) 
\end{eqnarray*}

But 
\begin{equation*} 
\max \{t,1\}\left( \int\nolimits_{t^{2}}^{\infty }\left( s^{\frac{1}{2}%
}f^{\ast \ast }\left( s\right) \right) ^{2}\frac{ds}{s}\right) ^{\frac{1}{2}%
}=\left\{ 
\begin{array}{c}
\left( \int\nolimits_{t^{2}}^{\infty }\left( s^{\frac{1}{2}}f^{\ast \ast
}\left( s\right) \right) ^{2}\frac{ds}{s}\right) ^{\frac{1}{2}}\text{ if }%
0<t<1 \\ 
t\left( \int\nolimits_{t^{2}}^{\infty }\left( s^{\frac{1}{2}}f^{\ast \ast
}\left( s\right) \right) ^{2}\frac{ds}{s}\right) ^{\frac{1}{2}}\text{ if }t>1%
\end{array}%
\right. 
\end{equation*}%
. Then, for $t<1$ it is 
\begin{eqnarray*} 
\left( \int\nolimits_{t^{2}}^{\infty }\left( s^{\frac{1}{2}}f^{\ast \ast
}\left( s\right) \right) ^{2}\frac{ds}{s}\right) ^{\frac{1}{2}} &\leq
&\lim\limits_{t\rightarrow 0^{+}}\left( \int\nolimits_{t^{2}}^{\infty
}\left( s^{\frac{1}{2}}f^{\ast \ast }\left( s\right) \right) ^{2}\frac{ds}{s}%
\right) ^{\frac{1}{2}} \\
&=&\left( \int\nolimits_{0}^{\infty }\left( f^{\ast \ast }\left( s\right)
\right) ^{2}ds\right) ^{\frac{1}{2}}=\left\Vert f\right\Vert _{L^{2}} 
\end{eqnarray*}%
.

And for $t>1$ we get $t\left( \int\nolimits_{t^{2}}^{\infty }\left( s^{-%
\frac{1}{2}}sf^{\ast \ast }\left( s\right) \right) ^{2}\frac{ds}{s}\right) ^{%
\frac{1}{2}}$; so, being $sf^{\ast \ast }\left( s\right)
=\int_{0}^{s}f^{\ast }\left( u\right) du$, for $s>t>1$ it is $sf^{\ast \ast
}\left( s\right) =\int_{0}^{s}f^{\ast }\left( u\right)
du=\int_{0}^{1}f^{\ast }\left( u\right) du=\left\Vert f\right\Vert
_{L^{1}\left( 0,1\right) }$ and then 
\begin{eqnarray*} 
t\left( \int\limits_{t^{2}}^{\infty }\left( s^{-\frac{1}{2}}\left\Vert
f\right\Vert _{L^{1}}\right) ^{2}\frac{ds}{s}\right) ^{\frac{1}{2}}
&=&\left\Vert f\right\Vert _{L^{1}}t\left( \int\limits_{t^{2}}^{\infty
}\left( s^{-\frac{1}{2}}\right) ^{2}\frac{ds}{s}\right) ^{\frac{1}{2}} \\
&=&\left\Vert f\right\Vert _{L^{1}}t\left( \int\limits_{t^{2}}^{\infty
}s^{-2}ds\right) ^{\frac{1}{2}}=\left\Vert f\right\Vert _{L^{1}}t\frac{1}{t}%
=\left\Vert f\right\Vert _{L^{1}} 
\end{eqnarray*}

That is 
\begin{eqnarray*} 
\left\Vert f\right\Vert _{\Delta _{0<\theta <1}\left( L^{\frac{1}{1-\frac{%
\theta }{2}}\blacktriangleleft }\left( 0,1\right) \right) } &\approx
&J\left( 1,f,L^{1},L^{2}\right) \\
&=&\max \{\left\Vert f\right\Vert _{L^{1}},\left\Vert f\right\Vert
_{L^{2}}\}=\left\Vert f\right\Vert _{L^{2}} 
\end{eqnarray*}%
and therefore the extrapolation result for the $\Delta $-functor: 
\begin{equation*} 
\Delta _{1<p<2}\left( \omega \left( p\right) L^{p\blacktriangleleft }\right)
=\Delta _{1<p<2}\left( \omega \left( p\right) \left( L^{1},L^{2}\right) _{%
\frac{2}{p^{\prime }},\infty ,K\left( \cdot ,\cdot ,L^{1},L^{2}\right)
}\right) 
\end{equation*}%
implies that 
\begin{eqnarray*} 
&&\inf\limits_{f\in \Phi }\left\Vert \chi _{\left( 0,1\right) }-f\right\Vert
_{\Delta _{1<p<2}\left( \omega \left( p\right) \left( L^{1},L^{2}\right) _{%
\frac{2}{p^{\prime }},\infty ,K\left( \cdot ,\cdot ,L^{1},L^{2}\right)
}\right) } \\
&\approx &\inf\limits_{f\in \Phi }\left\Vert \chi _{\left( 0,1\right)
}-f\right\Vert _{\Delta _{0<\theta <1}\left( L^{\frac{1}{1-\frac{\theta }{2}}%
\blacktriangleleft }\left( 0,1\right) \right) }\approx \inf\limits_{f\in
\Phi }\left\Vert \chi _{\left( 0,1\right) }-f\right\Vert _{L^{2}} 
\end{eqnarray*}%
,$\ $then for $M\left( \theta \right) \equiv 1$ the condition ii) $\chi
_{\left( 0,1\right) }\in \overline{\Phi }^{\Delta _{1<p<2}\left( \omega
\left( p\right) \left( L^{1},L^{2}\right) _{\frac{2}{p^{\prime }},\infty
,K\left( \cdot ,\cdot ,L^{1},L^{2}\right) }\right) }$ on the theorem \ref%
{strong Nyman-Beurling-Milman equivalence} becomes the same of the condition
ii) $\chi _{\left( 0,1\right) }\in \overline{\Phi }^{L^{2}}$ of the Nyman
-Beurling theorem \ref{N-B equivalence}.

\begin{remark}
The statement c) in our criterion \ref{N-B-with-extrapol} involve the $%
J-functional$ for the scale $\left\{ M\left( \theta \right) \overline{A}%
_{\theta ,\infty ,K}\right\} $ of interpolation spaces because 
\begin{equation*} 
J\left( t,\chi _{\left( 0,1\right) }-f,\left\{ M\left( \theta \right) 
\overline{A}_{\theta ,\infty ,K}\right\} \right) =\left\Vert \chi _{\left(
0,1\right) }-f\right\Vert _{\Delta _{0<\theta <1}\left( M\left( \theta
\right) \overline{A}_{\theta ,\infty ,K}\right) } 
\end{equation*}%
. It would be desirable to compute a formula for $J\left( t,f,\left\{
M\left( \theta \right) F_{\rho _{\theta }}\left( \overline{A}\right)
\right\} \right) $ for any tempered weight $M\left( \theta \right) $ and any
family of interpolation functors $\left\{ M\left( \theta \right) F_{\rho
_{\theta }}\left( \overline{A}\right) \right\} $ for quasi-concave functions 
$\rho _{\theta }$. This is an important open problem of Extrapolation Theory
(cf. \cite{AM} Problem 4.10.).

Of course, as we have already seen in the previous remark, we do know $J$
(up to equivalence) for the standard family $\left\{ L^{\frac{1}{1-\frac{%
\theta }{2}}\blacktriangleleft }\left( 0,1\right) \right\} _{\theta \in
\left( 0,1\right) }$.
\end{remark}

\section{The necessary condition and Reverse H\"{o}lder Classes}

Now we want to return to the question of the relationship between the
alleged invalidity of the Riemann hypothesis and the Reverse H\"{o}lder
weights in $\left( 0,1\right) $.

We will need a theorem due to Ivanov and Kalton (Ref. \cite{IK}) which is a
strengthening of a previous result of L\"{o}fstr\"{o}m (Cf. \cite{L}). An
even more general result is obtained in \cite{ACK} (see also \cite{AS})
which notation we adopt:

For $\left( Y_{0},Y_{1}\right) $ a regular couple, let 
\begin{equation*} 
\alpha _{0}\left( \psi \right) =\sup \left\{ \theta \in \left[ 0,1\right] :%
\frac{K\left( s,\psi ,Y_{0}^{\ast },Y_{1}^{\ast }\right) }{K\left( t,\psi
,Y_{0}^{\ast },Y_{1}^{\ast }\right) }\leq \gamma \left( \frac{s}{t}\right)
^{\theta }\text{ for all }0<s<t\leq 1\right\} 
\end{equation*}
for some $\gamma >0$ independent of $s$ and $t$; and 
\begin{equation*} 
\beta _{0}\left( \psi \right) =\inf \left\{ \theta \in \left[ 0,1\right] :%
\frac{K\left( s,\psi ,Y_{0}^{\ast },Y_{1}^{\ast }\right) }{K\left( t,\psi
,Y_{0}^{\ast },Y_{1}^{\ast }\right) }\geq \gamma \left( \frac{s}{t}\right)
^{\theta }\text{ for all }0<s<t\leq 1\right\} 
\end{equation*}
for some $\gamma >0$ independent of $s$ and $t$. It's clear that $K\left( \cdot ,\psi ,Y_{0}^{\ast },Y_{1}^{\ast }\right) $ is
increasing and $\frac{K\left( s,\psi ,Y_{0}^{\ast },Y_{1}^{\ast }\right) }{s}
$ is decreasing, and for $\theta <\alpha _{0}\left( \psi \right) $ we have
that $K\left( s,\psi ,Y_{0}^{\ast },Y_{1}^{\ast }\right) \cdot s^{-\theta }$
is almost increasing, which is the definition of index $i\left( \phi \right) 
$ for a quasi-concave function given in (Corval\'{a}n-Milman \cite{CM}
(4.3), see also Karapetyans-Samko \cite{KS})).

\begin{theorem}
(Ivanov-Kalton)

Suppose that $\left( X_{0},X_{1}\right) $ and $\left( Y_{0},Y_{1}\right) $
are couples of Banach spaces satisfying the following conditions:

i) the couple $\left( Y_{0},Y_{1}\right) $ is regular\footnote{%
This means that $Y_{0}\cap Y_{1}$ is dense in $Y_{0}$ and $Y_{1}$.}, $%
X_{0}\cap X_{1}$ is dense in $X_{1}$ and $X_{1}=Y_{1}$.

ii) $X_{0}$ is a closed subspace of codimension one of $Y_{0}$.

iii) $1\leq q<\infty $

Let $\psi \in Y_{0}^{\ast }$ be such that $X_{0}=\ker \psi $ then\ $\left(
X_{0},X_{1}\right) _{\theta ,q}$ is a closed subspace of codimension one of $%
\left( Y_{0},Y_{1}\right) _{\theta ,q}$ if and only if $\theta \notin \left[
\alpha _{0}\left( \psi \right) ,\beta _{0}\left( \psi \right) \right] $.
Moreover we have that $\left( X_{0},X_{1}\right) _{\theta ,q}$ is a closed
subspace of codimension one in $\left( Y_{0},Y_{1}\right) _{\theta ,q}$ for $%
\theta \in \left( 0,\alpha _{0}\left( \psi \right) \right) $, and the space $%
\left( X_{0},X_{1}\right) _{\theta ,q}$ coincides with $\left(
Y_{0},Y_{1}\right) _{\theta ,q}$ for $\theta \in \left( \beta _{0}\left(
\psi \right) ,1\right) $.
\end{theorem}

Now, let $Y_{1}=X_{1}=L^{1}\left( 0,1\right) $ and $Y_{0}=L^{2}\left(
0,1\right) $. If the Riemann Hypothesis is false then the closure of $\Phi $
in $L^{2}\left( 0,1\right) $ is a proper closed subspace of $L^{2}\left(
0,1\right) $ and then we can select $X_{0}$ a closed subspace of $%
L^{2}\left( 0,1\right) $ of codimension one such that $X_{0}\supseteq \Phi $%
. We know that $\Phi $ is dense in $L^{1}\left( 0,1\right) \footnote{%
This is a consequence of the existence of sequences in $\Phi $ that converge
to $\chi _{\left( 0,1\right) }$ in $L^{1}$ (cf. \cite{BDu4}) and the
equivalence of this fact with the density of $\Phi $ in $L^{1}$ (cf. \cite{B}%
).}$ and then $X_{1}=\overline{\Phi }^{L^{1}}$, and if $\Phi $ were dense in 
$L^{p}\left( 0,1\right) $ for all $1<p<2$ the Riemann Hypothesis would be
true, then if we assume that RH is false there is some $1<p<2$ such that $%
\left( X_{0},X_{1}\right) _{1-\theta ,\frac{2}{2-\theta }}$ doesn't coincide
with $\left( Y_{0},Y_{1}\right) _{1-\theta ,\frac{2}{2-\theta }}=\left(
Y_{1},Y_{0}\right) _{\theta ,\frac{2}{2-\theta }}=\left( L^{1}\left(
0,1\right) ,L^{2}\left( 0,1\right) \right) _{\theta ,\frac{2}{2-\theta }}=L^{%
\frac{2}{2-\theta }}\left( 0,1\right) =L^{p}\left( 0,1\right) $ for $p=\frac{%
2}{2-\theta }\in \left( 1,2\right) $ with $\theta \in \left( 0,1\right) $
and then $\alpha _{0}\left( \psi \right) >\theta >0$ for some $\psi \in
Y_{0}^{\ast }=\left( L^{2}\left( 0,1\right) \right) ^{\ast }=L^{2}\left(
0,1\right) $. Then, using that $\left( L^{1}\right) ^{\ast }=L^{\infty }$ we
have seen that for some $\psi \in L^{2}\left( 0,1\right) $ with $\ker \psi
=X_{0}\supseteq \overline{\Phi }^{L^{2}}$ we have that $i\left( K\left(
s,\psi ,L^{2}\left( 0,1\right) ,L^{\infty }\left( 0,1\right) \right) \right)
>0$, and by Theorem 3 and Remark 4 of Corval\'{a}n-Milman (\cite{CM}) this
means that $\psi \in RH\left( L^{2}\left( 0,1\right) ,L^{\infty }\left(
0,1\right) \right) $ the abstract H\"{o}lder class for the pair $\left(
L_{\left( 0,1\right) }^{2},L_{\left( 0,1\right) }^{\infty }\right) $ (see 
\cite{CM} for the definition of abstract H\"{o}lder classes) and $\ker \psi
\supseteq \overline{\Phi }^{L^{2}}$ implies that $\int_{0}^{1}\psi \left(
x\right) \left\{ \frac{1}{ax}\right\} dx=0$ for $a>0$. So we have proven:

\begin{proposition}
If the Riemann Hypothesis is false there is $\psi \in L^{2}\left( 0,1\right) 
$ such that $\psi \in RH\left( L^{2}\left( 0,1\right) ,L^{\infty }\left(
0,1\right) \right) $ and $\int_{0}^{1}\psi \left( x\right) \left\{ \frac{1}{%
ax}\right\} dx=0$ for all $a\geq 1$.
\end{proposition}

\begin{acknowledgement}
I would like to thank Mario Milman for the advice of properly clarify when
the necessary renormings are being considered to obtain functors which
characteristic functions are exactly $t^{\theta }\footnote{%
Among many other debts of gratitude.}$.
\end{acknowledgement}

\end{document}